\documentclass[sn-mathphys-num]{sn-jnl}

\usepackage{graphicx}%
\usepackage{multirow}%
\usepackage{amsmath,amssymb,amsfonts}%
\usepackage{amsthm}%
\usepackage[title]{appendix}%
\usepackage{xcolor}%
\usepackage{textcomp}%
\usepackage{manyfoot}%
\usepackage{booktabs}%
\usepackage{algorithm}%
\usepackage{algorithmicx}%
\usepackage{algpseudocode}%
\usepackage{listings}%
\usepackage[caption = false]{subfig}

\usepackage{color}

\usepackage{stmaryrd}
\usepackage[percent]{overpic}
\usepackage{array}
\usepackage{lmodern}

\newtheorem{theorem}{Theorem}
\newtheorem{lemma}{Lemma}

\graphicspath{{figures/}}

\begin{document}

\title[SGIGA2]{Stable quadratic generalized IsoGeometric analysis for elliptic interface problem}

\author[1]{\fnm{Yin} \sur{Song}}

\author[2]{\fnm{Wenkai} \sur{Hu}}

\author*[1]{\fnm{Xin} \sur{Li}}\email{lixustc@ustc.edu.cn}

\affil[1]{\orgdiv{School of Mathematical Science}, \orgname{USTC}, \orgaddress{\street{Jinzhai Road No.96}, \city{Hefei}, \postcode{230026}, \country{China}}}

\affil[2]{ \orgname{Jinhang Digital Technology Co. Ltd}, \orgaddress{\city{Beijing}, \postcode{100028}, \country{China}}}

\abstract{{Unfitted mesh formulations for interface problems generally adopt two distinct methodologies: (i) penalty-based approaches and (ii) explicit enrichment space techniques. While Stable Generalized Finite Element Method (SGFEM) has been rigorously established for one-dimensional and linear-element cases, the construction of optimal enrichment spaces preserving approximation-theoretic properties within isogeometric analysis (IGA) frameworks remains an open challenge.} In this paper, we introduce a stable
quadratic generalized isogeometric analysis (SGIGA2) for two-dimensional elliptic interface problems. The method is achieved through two
key ideas: a new quasi-interpolation for the function with $C^0$ continuous along interface and a new enrichment space with
controlled condition number for the stiffness matrix. We mathematically prove that the present method has optimal convergence rates
for elliptic interface problems and demonstrate its stability and robustness through numerical verification.}




\maketitle

\section{Introduction}\label{sec1}

{The key concept of isogeometric analysis (IGA) is the utilization of basis functions derived from Computer-Aided Design (CAD) models for the approximation of partial differential equations (PDEs).  This approach serves as a bridge between traditional engineering design and analysis methods.} In this context, shape functions extracted from CAD models typically serve as the initial approximation basis for solving PDEs. Accurate solutions can be obtained through global refinement ~\cite{HUGHES20054135} or local refinement~\cite{GIANNELLI2012485,WEI2022114494,BAZILEVS2010229,PATRIZI2020113230} if the solution to the PDEs is sufficiently smooth.
This approach will fail if the solution to the PDEs has reduced continuity, such as in the case of interface problems.

{Interface problems emerge from material property discontinuities across domains with heterogeneous physical characteristics. Strong discontinuity corresponds to $C^{-1}$ continuous solutions at the interface, whereas weak discontinuity describes $C^0$
continuous solutions.} Interface problems have attracted considerable attention in various fields, such as closed interfaces in solid
mechanics~\cite{MOES20033163}, moving interfaces for immiscible two-phase flows in fluid dynamics~\cite{SAUERLAND20113369,otmani2024accelerating},
and open interfaces in material analysis~\cite{RETHORE20075016,shakur2024isogeometric,gholampour2021global}.

{The numerical treatment of interface problems predominantly employs two distinct methodologies: conforming (fitted) mesh versus non-conforming (unfitted) mesh approaches.} The fitted mesh method requires the creation of a conforming mesh that aligns with the interface, which can be computationally expensive, particularly for dynamic interfaces that necessitate frequent mesh updates.  In contrast, unfitted mesh methods eliminate the need for such adjustments, offering potential computational efficiency benefits.  
This advantage arises because the unfitted mesh method only modifies the local basis functions influenced by the interface. The unfitted mesh method mainly includes two
approaches: the penalty method and the enrichment space method. The unfitted mesh method with penalty in general splits each basis function cross the interface into two basis functions and performs penalty term on these split basis functions, which is firstly introduced in~\cite{nitsche}.
The jump condition on the interface is enforced by penalty \cite{harari2010analysis,huang2017unfitted,chen2023arbitrarily,lehrenfeld2013analysis}.
{The primary challenge in penalty method implementation is the small cut cell problem:} the cut cells can be arbitrarily small and anisotropic,
which can make the stiffness matrix extremely ill-conditioned, especially for high-order finite element methods~\cite{burman2010fictitious,burman2012fictitious,huang2017unfitted,johansson2013high}.
The other unfitted mesh method is using enrichment space to approximate the solution with {reduced continuity} using the
generalized finite element method (GFEM) \cite{DENG2020112558,BABUSKA201758,ZHANG2019538,ZHANG2020112926,ZHU2020112475,GONG2024115540} or
generalized isogeometric analysis (GIGA)~\cite{8172107,ZHANG2022114053,TAN20151382}. GFEM and GIGA methods have been widely used in engineering problems~\cite{KIENDL20093902,GOMEZ20084333,Cazzani2016IsogeometricAO,THAI2013196,BabuOsborn_2003,fries2010extended,ABDELAZIZ20081141,iqbal2022enriched}.
{The unfitted mesh method with a penalty term splits some of the CAD basis functions into two distinct functions, altering the original CAD basis functions and violating the Iso-geometry principle in IGA. Therefore, the main motivation of this paper is to address the interface problem by generalizing the concept of GIGA, preserving all the desirable properties while maintaining the integrity of the CAD basis functions.}

For GFEM or GIGA, three main requirements are: (1) optimal convergence rates, (2) stability, where the conditioning should be stable
during h-refinement and p-refinement, and (3) Robustness, where the conditioning should be stable when the interfaces approach the
boundaries of some elements. {Generally, the GFEM method and GIGA method that meet the above three requirements are called
stable generalized FEM (SGFEM) and stable generalized IGA (SGIGA).} \cite{BABUSKA201291} proposed a stable generalized FEM method of
degree one with unfitted mesh to deal with the two-dimensional interface problem. However, {the similar idea} fails to meet
the three requirements to higher degrees. Lots of researches {try to improve condition number of GFEM by using orthogonalization}~\cite{AGATHOS20191051,schweitzer2011stable} or subtracting interpolation~\cite{BABUSKA201758,babu94,SAUERLAND201341,GUPTA201323,ZHANG2016476}.
The others try to improve the convergence by considering the blending elements~\cite{article,Cheng2009HigherorderXF} or imposing the boundary
conditions~\cite{HOU2020113135}. \cite{DURGARAO2019535,ZHANG2022114053,HU2024115792} achieve the goal
only when the interface is parallel to the parametric lines. To the best of the authors' knowledge, the other construction that satisfies the
three requirements for smooth interfaces is the stable GFEM of degree two (SGFEM2) in~\cite{ZHANG2020112889}. The paper proposes a stable GFEM of degree
two for the interface problem, where the finite element space is $C^0$ bi-quadratic piecewise polynomial space.

This paper focuses on the unfitted mesh method within the IGA framework, where the solution space is defined using traditional non-uniform B-splines. In contrast to the FEM space, the non-uniform B-spline space allows for a maximum continuity of  $C^{d-1}$  for B-splines of bi-degree  $d$. As continuity requirements increase, it becomes increasingly challenging to develop a method that ensures optimal convergence, stability, and robustness. To address these challenges, we {improve the idea in~\cite{ZHANG2020112889}} and introduce a stable quadratic generalized isogeometric analysis method (SGIGA2) for two-dimensional interface problems. The method is developed through the following key ideas: a new quasi-interpolation for the function with $C^0$ continuous along the interface and a new enrichment space with controlled condition number for the stiffness matrix. We also prove mathematically that the method yields optimal convergence rates for the interface problems with smooth interface and numerically verified that the method is stable and robust. {In the summary, the present paper has the following several main contributions compared with those
in~\cite{ZHANG2020112889},
\begin{itemize}
  \item SGFEM2 is a special case of SGIGA2. The basis functions in SGFEM2 are piecewise quadratic polynomials with $C^0$ continuity, whereas those in the present work are non-uniform quadratic B-splines. When multiple knots are used, the basis functions presented in this paper reduce to those in SGFEM2.
  \item In our approach, we propose a new quasi-interpolation method for functions with interfaces, as well as a new enrichment space, both of which differ from those in~\cite{ZHANG2020112889}.
  \item {SGFEM2 necessitates local principal component analysis (LPCA) for stiffness matrix conditioning control. The convergence rates are highly sensitive to the parameter in LPCA.} In contrast, the enrichment space in our approach is explicitly defined, eliminating the need for LPCA.
\end{itemize}}

The structure of the rest of this paper is as follows. Section \ref{sec2} describes the concept of non-uniform rational B-splines and
the quasi-interpolation. The model problem is introduced in Section \ref{sec3} with the introduction of FEM solver. In Section \ref{sec4},
we describe how to solve the interface problem with IGA. Section \ref{sec5} provides our new ideas for defining SGIGA2 and proving that the
method has optimal convergence rates. Section \ref{sec6} provides a performance comparison between IGA, SGIGA, SGIGA2 and SGFEM2.
The last section summarizes and prospects the entire paper with discussion and future work.

\section{Non-uniform rational B-splines}\label{sec2}

\begin{figure}[htbp]
    \centering
    \subfloat[basis functions in FEM]{     
        \includegraphics[width=0.45\linewidth]{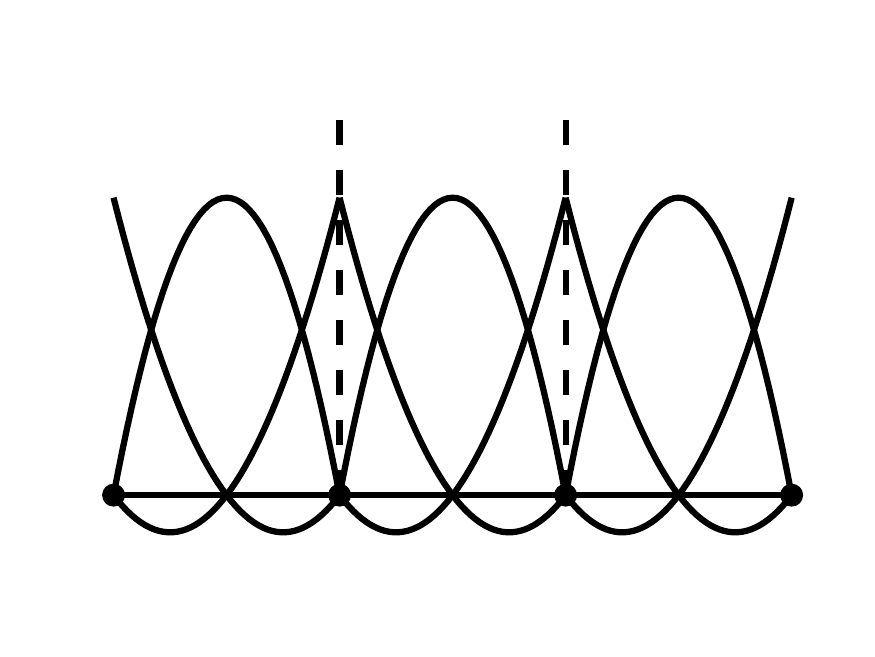}
        }\hfill
    \subfloat[basis functions in IGA]{
        \includegraphics[width=0.45\linewidth]{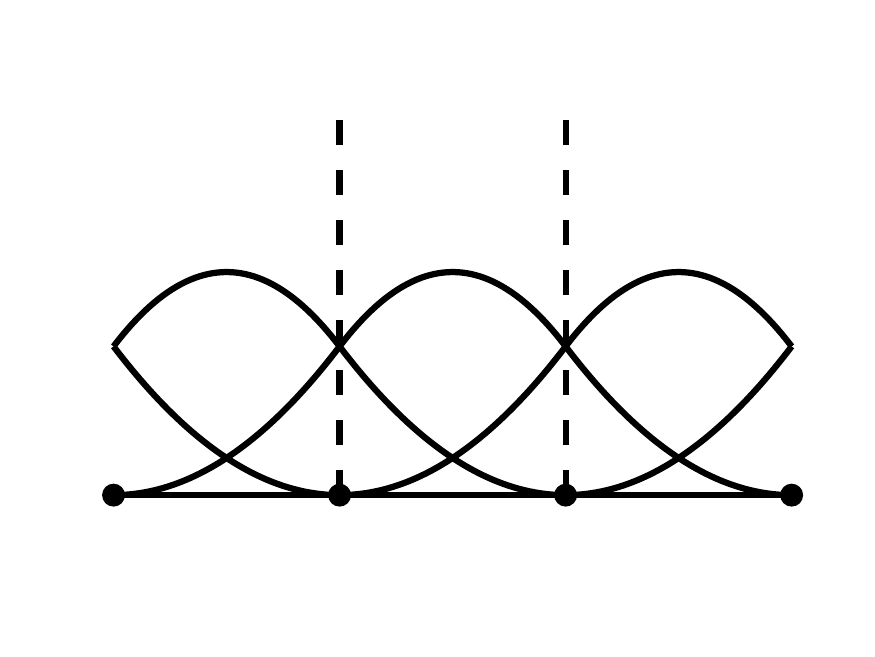}
        }
    \caption{Illustration of the quadratic basis functions in FEM (left) and IGA (right).}
    \label{FEM and B-spline basis function}
\end{figure}

{This section establishes the mathematical foundation for non-uniform rational B-Splines (NURBS) and their associated quasi-interpolation operators, which constitute the core components of our IGA framework.} B-spline basis function {$N^{d}_j(s)$} is defined in terms of a degree {$d$} and a set of knots $\{s_{j}\}$ with a recursion formula, where
{\begin{equation}
\begin{aligned}
& N_j^0(s)= \begin{cases}1, & \text { if } s_j \leq s < s_{j+1}, \\
0, & \text { otherwise, }\end{cases} \\
& N_j^d(s)=\frac{s-s_j}{s_{j+d}-s_j} N_j^{d-1}(s)+\frac{s_{j+d+1}-s}{s_{j+d+1}-s_{j+1}} N_{j+1}^{d-1}(s).
\end{aligned}
\end{equation}}
{The B-spline basis function $N^{d}_j(s)$ requires $d+2$ knots $s_{j}, \dots, s_{j+d+1}$.}
{Fig~\ref{FEM and B-spline basis function} depicts quadratic Lagrangian finite element basis functions versus $C^1$ continuous quadratic B-spline basis functions, with dashed vertical lines indicating element boundaries defined by nodal points (FEM) or knot spans (B-splines).}

With the definition of B-spline basis functions, we can define the NURBS curve and the surface as follows. Given a
degree $d$ and a knot sequence $\{s_{0}, s_{1}, \dots, s_{m+d-1}, s_{m+d}\}$, where $s_{0} = \dots = s_{d}$ and
$s_{m} = \dots = s_{m+d}$, the knot sequence can
define the $m$ degree $d$ B-spline basis functions $N_{j}^{d}(s)$, $j = 0, \dots, m-1$, then a NURBS curve is defined as
{\begin{equation}
C(s) = { \frac{\sum_{i = 0 }^{m-1} c_{i} \omega_{i} N_{i}^{d}(s)}{\sum_{i=0}^{m-1} \omega_{i} N_{i}^{d}(s) }}, \label{nncurve}
\end{equation}}
where {$c_{i}$} are control points, $\omega_{i} \geq 0$ are weights, $N_{i}^{d}(s)$ are B-spline basis function.

A NURBS surface is defined in terms of tensor-product setting. Given two degrees $p$, $q$, two knot sequences
$\{s_{0}, s_{1}, \dots, s_{m+p-1}, s_{m+p}\}$ and $\{t_{0}, t_{1}, \dots, t_{n+q-1}, t_{n+q}\}$ can
define a set of B-spline basis functions $N_{i}^{p}(s)$ and $N_{j}^{q}(t)$. A NURBS surface is defined as
{\begin{equation}
C(s, t) = { \frac{\sum_{i = 0 }^{m-1}\sum_{j = 0 }^{n-1} c_{i, j} \omega_{i, j} N_{i}^{p}(s)N_{j}^{q}(t)}
{\sum_{i=0}^{m-1}\sum_{j=0}^{n-1} \omega_{i, j} N_{i}^{p}(s)N_{j}^{q}(t)}}, \label{nnsurf}
\end{equation}}
where {$c_{i,j}$} are control points, $\omega_{i, j} \geq 0$ are weights, $N_{i}^{p}(s)$ and $N_{j}^{q}(t)$
are B-spline basis function. {For simplicity, $N_{i}(s)$ denotes $N_{i}^{2}(s)$ in the following. For quadratic B-splines, the B-spline basis functions correspond to the elements defined by the knot sequence.}

{We present a systematic formulation of quadratic B-spline quasi-interpolation operators, building upon the foundational works~\cite{xinli_quasi,barrera2020non}. The surface construction constitutes a straightforward extension via tensor product representations.} The following focuses on the detailed equations for the quasi-interpolation of a B-spline curve. Given knot sequence $\{s_{0}, s_{1}, \dots, s_{m+1}, s_{m+2}\}$,
suppose the quadratic B-spline basis functions are $N_{i}(s)$, $i = 0, \dots, m-1$, then the B-spline quasi-interpolation $I_{b}(\cdot)$
is a functional to project a function $f(s)\in L^{2}$ into B-spline space, i.e.,
\begin{equation}
I_{b}(f(s)) = \sum_{i = 0}^{ m - 1} \mu_i(f(s)) N_i(s),
\end{equation}
where the coefficients $\mu_j(f(s))$ are linear combination of the value of $f(s)$ at three points
$\tau_{j}^{k} = \frac{s_{j+k}+s_{j+k+1}}{2}$, $k = 0, 1, 2$, i.e.,
\begin{equation}
\mu_{j, 3}(f) =\alpha_{j, 0} f(\tau_{j}^{0})+\alpha_{j, 1} f(\tau_{j}^{1})+\alpha_{j, 2}f(\tau_{j}^{2}),
\end{equation}
where
\begin{equation}
\alpha_{j, 0} =-\frac{\tau_j^2 (1 - \tau_{j+1})}{1 + \tau_j - \tau_{j+1}}, \quad
\alpha_{j, 1} =1+\tau_j (1-\tau_{j+1}), \quad
\alpha_{j, 2} =-\frac{\left(1 - \tau_{j+1}\right)^2 \tau_j}{1 + \tau_j - \tau_{j+1}}, \quad
\tau_j:=\frac{s_{j+2}-s_{j+1}}{s_{j+2}-s_{j}}.
\end{equation}

\begin{figure}[htbp]
    \centering
    \subfloat[The Fig shows the results of quasi-interpolation for a smooth function, as well as the three values of $f$ (green points:  $f(\tau_{j}^{0})$,$f(\tau_{j}^{1})$,$f(\tau_{j}^{2})$) needed to calculate a B-spline coefficient (the black dot represents a B-spline control point and the arrows represent the support of this B-spline).]{     
        \includegraphics[width=0.45\linewidth]{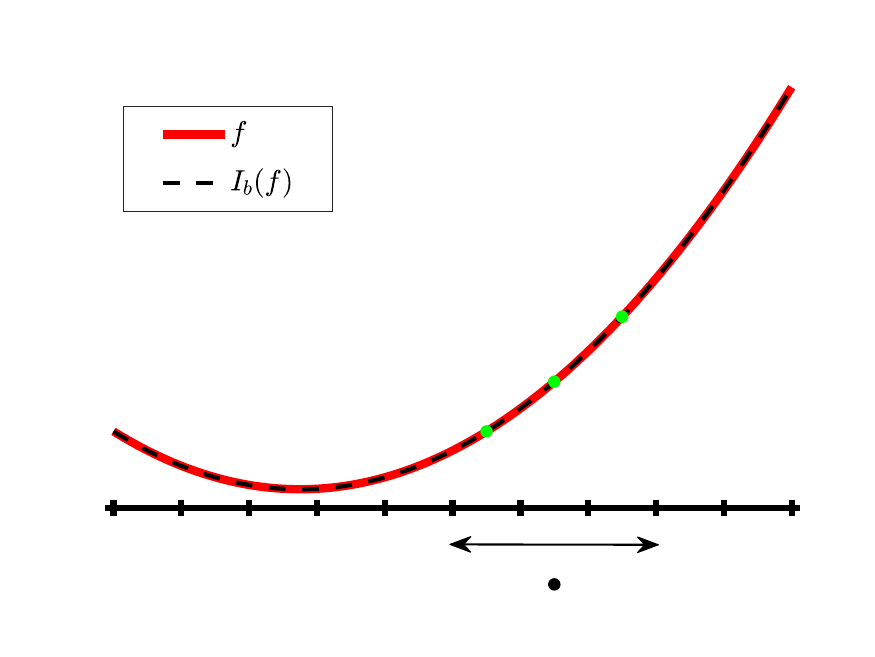}
        }\hfill
    \subfloat[The Fig shows the results of quasi-interpolation for a function that is $C^0$ continuous at the interface, where the error theorem does not apply (indicated by the red elements).]{
        \includegraphics[width=0.45\linewidth]{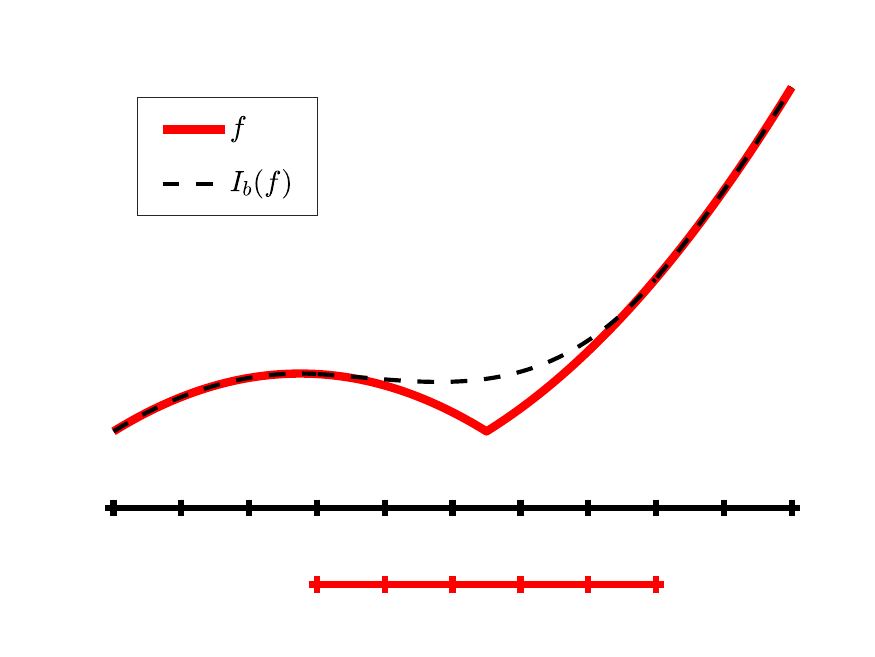}
        }
    \caption{Illustration of the quasi-interpolation for quadratic B-splines.}
    \label{quasi_1}
\end{figure}

\begin{theorem}\label{quasierror}
	There exists a constant $0<C<1$ such that for all $u \in W^{3, \infty}(\Omega)$  with step size $h$, we have
	\begin{equation}
	\left\|u-I_b u\right\|_{L^{\infty}(\Omega)} \leqslant C h^3\left\|u\right\|_{W^{3, \infty}(\Omega)}
	\end{equation}
\end{theorem}
{The aforementioned quasi-interpolation scheme delivers optimal approximation rates for functions in $W^{3, \infty}(\Omega)$, as shown in \cite{barrera2020non}.} The above quasi-interpolation works well for smooth function. However, for the function with interface, we find that the above theorem cannot be applied for five elements (red elements). Fig \ref{quasi_1} shows the results of quasi-interpolation for a quadratic polynomial and a piecewise quadratic polynomial with $C^0$ continuous respectively.

\section{Interface problem and solving framework}\label{sec3}

In this section, we present the interface problem and the corresponding solver using the GFEM or GIGA.

Let $\Omega$ be a bounded and simply connected domain with piecewise smooth boundary ${\partial \Omega}$. For a domain $\Omega$ in $\mathbb{R}^2$, an integer $m$ and $1 \leq q \leq \infty$, we denote sobolev spaces as $W^{m, q}(\Omega)$ with norm ${\|\cdot\|_{W^{m, q}(\Omega)}}$ and semi-norm ${|\cdot|_{W^{m, q}(\Omega)}}$. The space $W^{m, q}(\Omega)$ will be represented by $H^m(\Omega)$ for $q=2$ and $L^q(\Omega)$ when $m=0$, respectively. We assume that $\Omega$ has a smooth interface $\Gamma$ that divides
$\Omega$ into two domains $\Omega_0$ and $\Omega_1$ such that $\bar{\Omega}=\bar{\Omega}_0 \cup \bar{\Omega}_1,
\Omega_0 \cap \Omega_1=\emptyset$, and $\Gamma=\bar{\Omega}_0 \cap \bar{\Omega}_1$.
We consider the following second-order elliptic problem on $ \Omega $ with Neumann  boundary conditions,
\begin{equation}\label{problem}
	 \begin{aligned}
	&-\nabla \cdot(a \nabla u)=f, \text { in } \Omega \text {; }\\
	&a \nabla u \cdot \vec{n}_b=g \quad \text {, on } \partial \Omega \text {, }\\
	& {\llbracket u \rrbracket_{\Gamma}=0 \text {, on } \Gamma \text {, }} \\
	& { \llbracket a \nabla u \cdot \boldsymbol{n} \rrbracket }_{\Gamma}=q \text {, on } \Gamma \text {, } \\
	\end{aligned}
\end{equation}
where $\llbracket \cdot \rrbracket_{\Gamma}$ denotes the jump of the underlying function across the interface $\Gamma$, defined as
\begin{equation}
\begin{aligned}
\llbracket v \rrbracket _{\Gamma}  &:=\left.v\right|_{\Omega_0}-\left.v\right|_{\Omega_1}, \quad \\
\quad \llbracket \nabla v \cdot \boldsymbol{n} \rrbracket&:=\left.\nabla v\right|_{\Omega_0} \cdot \boldsymbol{n}_{\Omega_0}+\left.\nabla v\right|_{\Omega_1} \cdot \boldsymbol{n}_{\Omega_1} \text {, }
\end{aligned}
\end{equation}
The given functions $f, g$, and $q$ satisfy a compatibility condition:
\begin{equation}
\int_{\Omega} f d P+\int_{\partial \Omega} g d s+\int_{\Gamma} q d s=0 .
\end{equation}
{The piecewise-constant function $a$ given by
\begin{equation}
    a =  \begin{cases}a_0, & (x,y)\in \Omega_0, \\ a_1, & (x,y)\in \Omega_1 .\end{cases}.
\end{equation}
The notation $\vec{n}_b$ denote the unit outward normal to the boundary $\partial \Omega$.} The variational formulation of interface problem can be defined as that
for $f \in L^2(\Omega)$, find $u \in \{H^1(\Omega)| [v]_{\Gamma}=0 \text { on } \Gamma\}$ such that
\begin{equation}
B(u, v)=L(v), \quad \forall v \in \{H^1(\Omega)| [v]_{\Gamma}=0 \text { on } \Gamma\},
\end{equation}
where
\begin{equation}
\begin{aligned}
B(u, v)&:=\int_{\Omega}a \nabla u \cdot \nabla v d{\Omega}, \\
L(v)&:=\int_{\Omega} f v
d {\Omega}+\int_{\partial \Omega} g v d s+\int_{\Gamma} q v d s.
\end{aligned}
\end{equation}

{The GFEM and GIGA employ enrichment spaces to locally resolve solution discontinuities, augmenting standard approximation spaces through direct summation} $V^h=V_{IGA/FEM}^h \oplus V_{ENR}^h$, where $V_{IGA/FEM}^h $ is the standard FEM or
IGA space and $V_{ENR}^h$ is the enrichment space, i.e., the approximation solution
\begin{equation}
u^h(s, t) = \sum_{i\in \mathcal{J}^{h}} c_{i} \phi_i(s, t)+\sum_{j\in \mathcal{J}^{h}_{E}} d_j \psi_j(s, t),
\end{equation}
where
$\phi_{i}(s, t)$ and $\psi_j(s, t)$ are basis functions for $V_{IGA/FEM}^h$ and $V_{ENR}^h$. $\mathcal{J}^{h}$ and $\mathcal{J}^{h}_{E}$ are the indices sets for the original space and the enrichment space respectively.

Let {${ \boldsymbol{c} }$} and {${ \boldsymbol{d} }$} be the vectors of the coefficients
$c_i$ and $d_j$ and they can be determined by the system of linear equations
 \begin{equation}
 \left[\begin{array}{ll}
 \boldsymbol{K}_{O O} & \boldsymbol{K}_{O E} \\
 \boldsymbol{K}_{E O} & \boldsymbol{K}_{E E}
 \end{array}\right]\left[\begin{array}{l}
{\boldsymbol{c}}  \\
{\boldsymbol{d}}
 \end{array}\right]=\left[\begin{array}{l}
 {\boldsymbol{f}_O} \\
 {\boldsymbol{f}_E}
 \end{array}\right],
\end{equation}
 where the block matrices are defined with entries
 \begin{equation}
 \begin{aligned}
 & \left(\boldsymbol{K}_{O O}\right)_{j k}=B\left(\phi_k, \phi_j\right), \left(\boldsymbol{K}_{E E}\right)_{j k}=B\left(\psi_k, \psi_j\right), \\
 & \left(\boldsymbol{K}_{O E}\right)_{j k}=\left(\boldsymbol{K}_{E O}\right)_{k j}=B\left(\psi_k, \phi_j\right),
 \end{aligned}
\end{equation}
and the load vector are defined as
\begin{equation}
({\boldsymbol{f}_O})_{j}=L\left(f, \phi_j\right), \quad ({\boldsymbol{f}_E})_{k}=L\left(f, \psi_k\right).
\end{equation}

{Since the stiffness matrix is positive-definite, the generalized matrix system guarantees a unique solution. In general, we will simplify the above equation as
\begin{equation}
\label{souce equation}
\boldsymbol{K} \boldsymbol{U}= \boldsymbol{F}
\end{equation}}

The numerical stability of a linear system is in generally measured by Scaled Condition Number (SCN) of $\hat{\boldsymbol{K}}$, where
\begin{equation}
\begin{gathered}
\hat{\boldsymbol{K}}=\boldsymbol{D} \boldsymbol{K} \boldsymbol{D}, \quad \hat{\boldsymbol{F}}=\boldsymbol{D} \boldsymbol{F}, \quad \hat{\boldsymbol{U}}=\boldsymbol{D}^{-1} \boldsymbol{U},
\end{gathered}
\end{equation}
and
\begin{equation}
\begin{aligned}
\boldsymbol{D}=\left[\begin{array}{cc}
\boldsymbol{D}_{11} & 0 \\
0 & \boldsymbol{D}_{22}
\end{array}\right], \quad\left(\boldsymbol{D}_{11}\right)_{i i}=\left(\boldsymbol{K}_{OO}\right)_{i i}^{-1 / 2}, \quad \left(\boldsymbol{D}_{22}\right)_{i i}=\left(\boldsymbol{K}_{EE}\right)_{i i}^{-1/2}.
\end{aligned}
\end{equation}
{This step is essential to mitigate the effects of rounding errors arising from numerical differences between the two spaces. See \cite{BABUSKA201291} for more details on scaling matrices. We solve a linear system equivalent to \eqref{souce equation} to avoid the possible loss of precision, which is
\begin{equation}
\hat{\boldsymbol{K}} \hat{\boldsymbol{U}}= \hat{\boldsymbol{F}}.
\end{equation}
}

\subsection{Interface problem using finite element method}
\begin{figure}[htbp]
    \centering
    \subfloat[The Fig shows $\mathcal{J}^{f}$ (yellow elements),
$\mathcal{J}^{f}_{1, +}$(blue and yellow elements) and $\mathcal{J}^{f}_{1, -}$(green and yellow elements)]{    
        \includegraphics[width=0.45\linewidth]{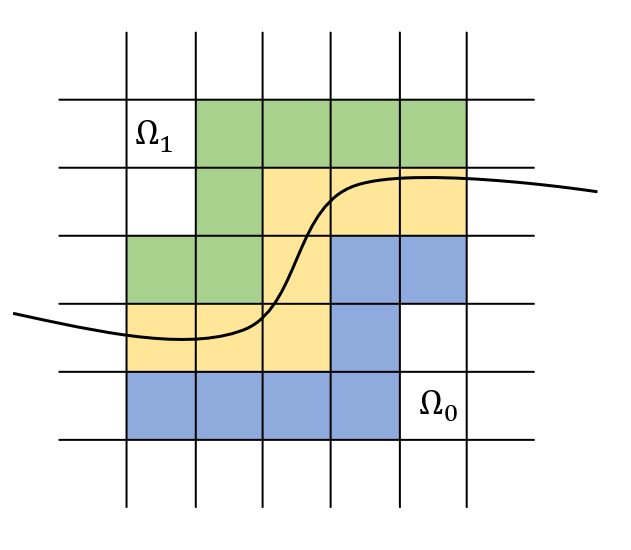}
        }\hfill
    \subfloat[The Fig shows $\mathcal{J}^{v}$ (red vertices) and $\mathcal{J}^{v}_{2}$ (green and red vertices)]{
        \includegraphics[width=0.45\linewidth]{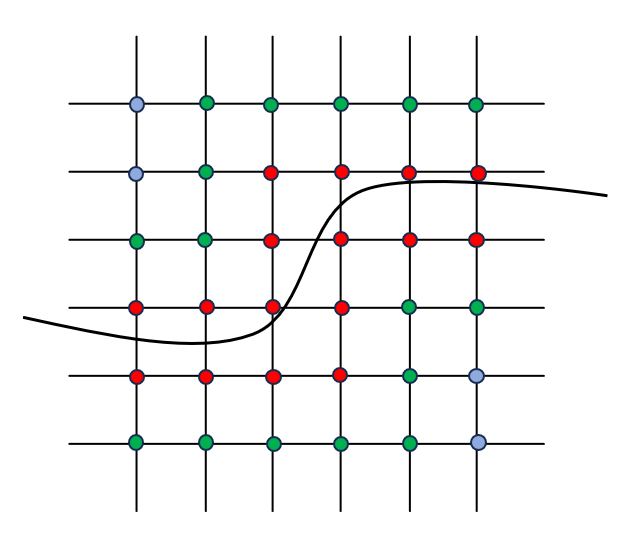}
        }
    \caption{Illustration of the notations for the construction of enrichment space.}
   \label{severalelements}
\end{figure}

{As demonstrated, resolving interface problems is fundamentally dependent on the appropriate construction of the finite space $ V_{IGA/FEM}^h $ and enrichment space $ V_{ENR}^h $.}
{
To systematically characterize interface-dependent topological components, we establish the following notational framework guided by the geometric configuration in Fig~\ref{severalelements}:}
{
\begin{itemize}
\item $\mathcal{J}^{f}_{k, \pm}$: $\mathcal{J}^{f}_{0}$ denotes indices of the elements intersect the interface, and $\mathcal{J}^{f}_{k}$
denotes the indices of all the elements in $\mathcal{J}^{f}_{k-1}$ and the elements which have at least one common vertices with the elements in
$\mathcal{J}^{f}_{k-1}$. And $\mathcal{J}^{f}_{k, +}$ denotes the indices of elements in $\mathcal{J}^{f}_{k}$ that have non-empty
intersections with $\Omega_0$ while $\mathcal{J}^{f}_{k, -}$ denotes the indices of elements in $\mathcal{J}^{f}_{k}$ that have non-empty
intersections with $\Omega_1$.
\item $\mathcal{J}^{v}_{k}$: $\mathcal{J}^{v}_{k}$ denotes the indices of all vertices in $\mathcal{J}^{f}_{k}$.
\end{itemize}
}

Now we are ready to review the existing construction of these two spaces using GFEM.
\begin{itemize}
  \item GFEM of degree one: $\phi_i$ is the hat function for the $j$-th vertex of the mesh and
  \begin{equation}V_{ENR}^h = \{\phi_{j} d_{\Gamma}(x,y):j\in\mathcal{J}^{v}\},\end{equation}
  where $d_{\Gamma}(x,y)$ is the distance function of a point to the interface $\Gamma$ for the weak discontinuous interface problems~\cite{haasemann2011development,kastner2013higher,BABUSKA201758,nme2914,ZHANG2019538,MOES20033163,BABUSKA201291}. However, the construction cannot yield optimal convergence rates and is not robust or stable. To handle this, SGFEM is introduced with the same space $V_{FEM}^{h}$ and different enrichment space,
  \begin{equation}V_{ENR}^h = \{\phi_{j} (d_{\Gamma}(x,y) - I_{h}(d_{\Gamma})(x,y)):j\in\mathcal{J}^{v}\},\end{equation}
  where $I_{h}(f)$ is the interpolation of function $f$ by $\phi_{i}$. {Through a carefully designed modification of the distance function, SGFEM of degree one is robust, stable and has optimal convergence rates.}
  \item Higher order GFEM: {While a straightforward extension from linear to higher-order elements appears theoretically plausible, significant technical challenges emerge in practical implementations.} For degree $d$, $V_{FEM}^{h}$ is defined as the piecewise degree $d$ polynomial spline space of
  $C^0$ continuous, and suppose $\phi_j$ is the $j$-th basis function of the space. Similarly, we can define
  $V_{ENR}^h = \{\phi_{j} d_{\Gamma}(x,y):j\in\mathcal{J}^{v}\}$ or $V_{ENR}^h = \{\phi_{j} (d_{\Gamma}(x,y) - I_{h}(d_{\Gamma})(x,y)):j\in\mathcal{J}^{v}\}$. However, both constructions
  do not have optimal convergence rates. The other ideas include adding more enrichment functions, for example,
  \begin{equation}V_{ENR}^h=\left\{\phi_{j}d_{\Gamma}(x, y)\left\{1,x,y,x^2,xy,y^2\right\}:j\in\mathcal{J}^{v}\right\}.\end{equation}
{In other words, there are actually six enrichment functions in higher order GFEM, which are $d_{\Gamma}(x, y)$, $xd_{\Gamma}(x, y)$, $yd_{\Gamma}(x, y)$, $x^2d_{\Gamma}(x, y)$, $xyd_{\Gamma}(x, y)$, $y^2d_{\Gamma}(x, y)$.}
  However, none of them can achieve the purposes of optimal convergence, robust and stable.
  \item Higher order Corrected GFEM: \cite{fries2008corrected} proposed a corrected GFEM method to improve convergence rates, which defines a modified enrichment function $d_{\Gamma}^{\text{mod}}(x,y) = d_{\Gamma}(x,y)r(x,y)$, where $r(x,y)$ is the ramp function.
  {There are several different methods to define the ramp functions in~\cite{fries2008corrected}, where the typical $r(x,y)$ is the sum of those basis functions that are not zero in the interface elements. However, all methods can improve the
  convergence rates, but none of them can achieve the optimal convergence.}
  \item SGFEM2: \cite{ZHANG2020112889} proposes a stable GFEM of degree two for interface problem, which yields optimal convergence. The
  enrichment space is
 \begin{equation}
 V_{ENR}^h=\left\{\varphi_j[I-I_h]\left(d_{\Gamma}(x,y)\left\{1,x,y\right\}\right) : j\in\mathcal{J}^{\text{v}} \right\},
 \end{equation}
where $I_h$ is piecewise FE of degree $2$ interpolation based on FEM basis function. And $\varphi_j$ is the Hermite function that satisfies the {partition of unity} properties, which is also defined by the mesh nodes. {In fact, there are three enrichment functions, which are $[I-I_h]\left(d_{\Gamma}(x,y)\right),[I-I_h]\left(d_{\Gamma}(x,y)x\right),[I-I_h]\left(d_{\Gamma}(x,y)y\right)$.}
\end{itemize}

{In the table \ref{different GFEM1}, we make a brief summary of several methods that appear in the literature.}
\begin{table*}[htbp]
    \centering
    \begin{tabular}{ccccccc}
   \toprule
   Prop. & degree & enrichment & Optimal convergence  & Stability & Robustness   \\
   \midrule
   FEM & $\geq 1$ & $\times$ & $\times$ & $\checkmark$ & $\checkmark$  \\
   GFEM & $1$ & $d_{\Gamma}$ & $\times$  & $\times$ & $\times$   \\
   SGFEM & $1$ & $[I-I_h](d_{\Gamma})$ & $\checkmark$ & $\checkmark$ & $\checkmark$   \\
Higher order GFEM & $\geq 2$ & $d_{\Gamma}$ & $\times$  &  $\times$ & $\times$ \\
Higher order Cor. GFEM & $\geq 2$& $d_{\Gamma}^{\text{mod}}$ & $\times$  &  $\times$ & $\times$\\
SGFEM2& $2$ & $[I-I_h]\left(d_{\Gamma}\left\{1,x,y\right\}\right)$ & $\checkmark$ & $\checkmark$ & $\checkmark$\\
   \bottomrule
\end{tabular}
    \caption{Comparison of GFEM methodologies.}
    \label{different GFEM1}
\end{table*}

\section{Interface problem using isogeometric analysis}\label{sec4}

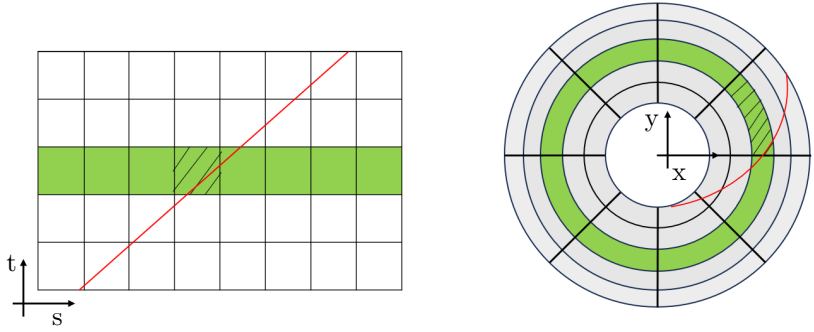
\begin{figure}[htbp]
\centering
\begin{overpic}
    [width=0.75\linewidth]{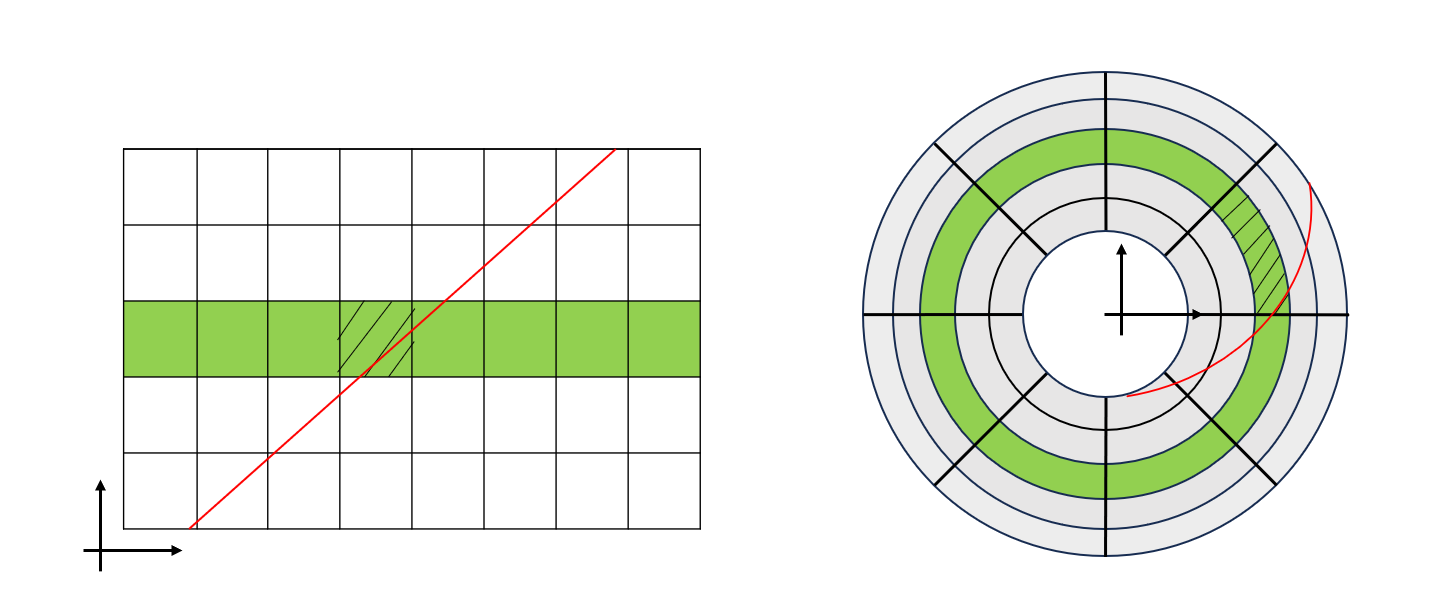}
     \put(10,2){s}
     \put(5,8){t}

     \put(78,18){x}
     \put(75,24){y}

 \end{overpic}
    \caption{Mapping between parameter domain and physical domain.}
    \label{mesh}
\end{figure}

{The fundamental premise of IGA lies in employing geometry-representing basis functions for PDE solution approximation.} As shown in Fig~\ref{mesh}, suppose the domain is presented using B-spline
basis functions as following,
\begin{equation}
    C(s, t)=\sum N_i^{d}(s, t) \cdot c_i(s, t),\quad {\Omega^p} \rightarrow {\Omega}
\end{equation}
where $N_i^{d}(s, t)$ are B-spline basis functions, $c_i$ are control points {and ${\Omega^p}$ is parameter domain}.

{The governing PDE solution procedure in IGA maintains conceptual alignment with conventional FEM frameworks.} In the IGA framework, the approximation space is typically defined as {$\{N_{i}^{d}\}$}. {The fundamental distinction between GIGA and GFEM lies in the enhanced continuity: GIGA employs $C^{d-1}$ continuous bidegree $d$ B-splines, whereas GFEM utilizes $C^{0}$ continuous piecewise polynomials. Furthermore,} GIGA is capable of handling curved domains while preserving the parameterization during the refinement process, as shown in the example~\ref{sec71}. {In the following, GIGA is introduced in~\cite{8172107} and
others are directly generalized from the previous section in the present paper.}
\label{GIGA}
\begin{itemize}
  \item GIGA: In the two-dimensional case, \cite{8172107} proposed GIGA using a distance function as enrichment to mimic the local behavior of the solution to the underlying variational problem, the approximation space is obtained by
      \begin{equation}V_{ENR}^h=\left\{d_{\Gamma}(s, t)  N_j(s, t) :j\in\mathcal{J}^{f}_1\right\}\end{equation}
      where the $N_j(s, t)$ are B-spline basis functions whose support contains the interface elements. For quadratic B-splines, each element corresponds to a B-spline basis function. So we use the previous notation $\mathcal{J}^{f}_1$ to represent the set of enriched B-splines.
      \item SGIGA: We can also borrow the idea from SGFEM to define the enrichment space as
      \begin{equation}
      \label{SGIGA}
      V_{ENR}^h = \{N_j(s, t) (d_{\Gamma} - I_{b}(d_{\Gamma})):j\in\mathcal{J}^{f}_1\},
      \end{equation}
     where the $I_b$ is quasi-interpolation.
  \item Corrected GIGA: we also try to generalize the idea of~\cite{fries2008corrected} to GIGA by defining the enrichment function as
  $d_{\Gamma}^{\text{mod}}(s, t) = d_{\Gamma}(s, t)r(s, t)$, where $r(s, t)$ is the sum of those $N_j(s, t)$ whose support contains the
  interface elements. However, unlike the GFEM case, corrected GIGA cannot improve convergence rates, as illustrated in Fig~\ref{the error of straight line interface}.
  \item SGIGA-multi: We also generalize the idea of \cite{ZHANG2020112889} to GIGA by the quasi-interpolation defined in section~\ref{sec2} and define
  \begin{equation}
  V_{ENR}^h=\left\{N_j(s, t)\left(\left(I - I_{b}\right)\left\{d_{\Gamma}, sd_{\Gamma}, td_{\Gamma}\right\}\right):j\in\mathcal{J}^{f}_1\right\},
 \end{equation}
  where $N_j(s, t)$ are B-spline basis functions whose support contains the interface elements.
\end{itemize}

The following example will verify the above conclusions.
\subsection{{An example of a straight line interface}}
The domain $ \Omega $ is chosen by $ \Omega := [0,1] \times [0,1] $ , which is divided into two parts by the straight line interface. The straight line has an equation $ y = \tan(\frac{\pi}{6})(x-1-\frac{1}{\pi}) + 1$. The manufactured solution of \eqref{problem} is as follows:
\begin{equation}
\begin{aligned}
 \begin{cases} u_0 = (r^{\frac{4}{3}}\cos(\frac{4}{3}(\theta +\pi - \frac{\pi}{6})) + \frac{a_0}{a_1}r^{\frac{4}{3}}\sin(\frac{4}{3}(\theta +\pi- \frac{\pi}{6})) \\+ \sin{xy} \text { if } y > \tan(\frac{\pi}{6})(x-1-\frac{1}{\pi}) + 1 (\Omega_0) \\
u_1 = (r^{\frac{4}{3}}\cos(\frac{4}{3}(\theta +\pi - \frac{\pi}{6})) + r^{\frac{4}{3}}\sin(\frac{4}{3}(\theta +\pi- \frac{\pi}{6}))\\
+ \sin{xy},  \text { if } y < \tan(\frac{\pi}{6})(x-1-\frac{1}{\pi}) + 1 (\Omega_1)
\end{cases}
\end{aligned}
\end{equation}
where the centre of polar coordinates $ (r,\theta) $ is $ (1+\frac{1}{\pi},1) $. It is easy to verify that this solution $ u $ is smooth on $ \Omega_0 $, $ \Omega_1 $ and $ C^0 $ at the interface.  {The domain is uniformly discretized into $N \times N$ elements with $ N = 5, 10, 20, \dots , 160 $}.
\begin{figure}[htbp]
    \centering
    \subfloat[The $H^1$ errors]{     
        \includegraphics[width=0.45\linewidth]{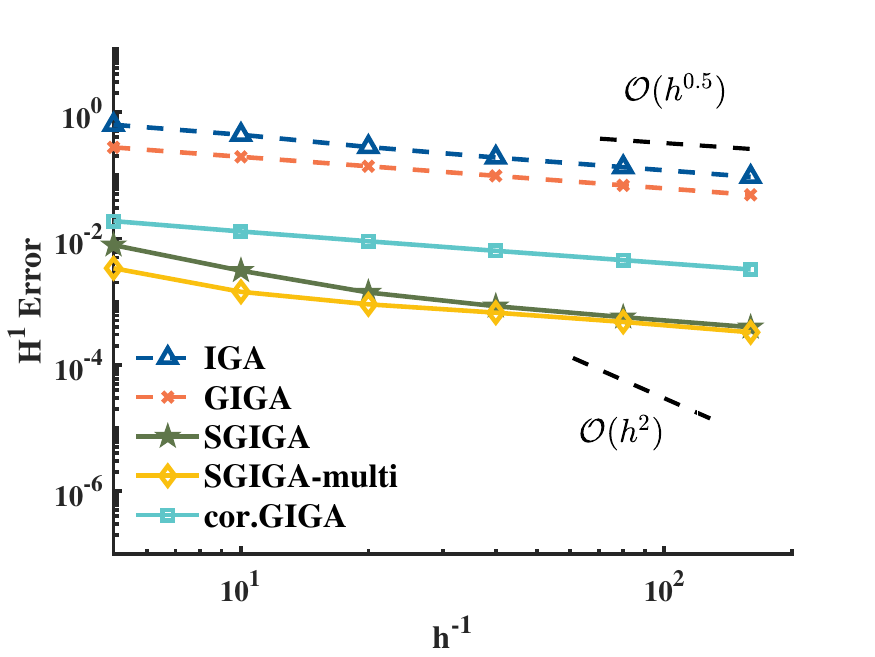}
        }\hfill
    \subfloat[The SCNs result]{
        \includegraphics[width=0.45\linewidth]{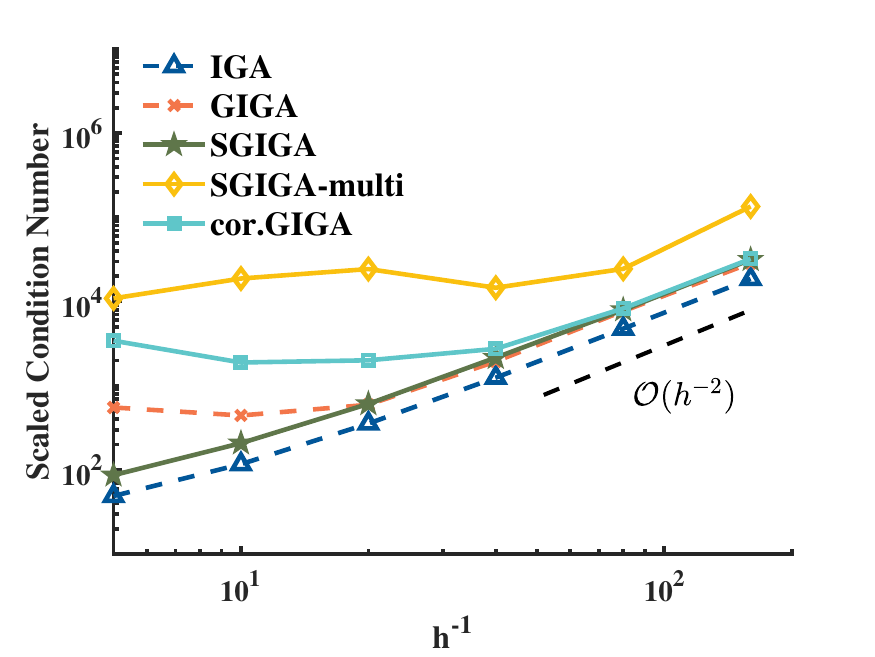}
        }
    \caption{The $H^1$ errors and the SCNs with respect to $h$ (from $\frac{1}{5}$ to
$\frac{1}{160}$) for a straight line interface ($ a_0 = 20, a_1 = 1$).}
    \label{the error of straight line interface}
\end{figure}

Fig~\ref{the error of straight line interface} present $ H^1 $ seminorm errors and SCNs of several methods for a straight line interface problem. It is clearly observed that the several methods do not {yield} optimal convergence. Compared with IGA, the error of other methods is improved, but the convergence speed is not significantly improved. The SCNs of the IGA, and SGIGA increase with the same order $ O(h^{-2}) $, which is shown in Fig~\ref{the error of straight line interface} while other several methods are not stable.
The table \ref{differentGIGA1} shows the result of several GIGA methods.
\begin{table*}[htbp]
    \centering
    \begin{tabular}{ccccccc}
   \toprule
   Prop. & enrichment & Optimal convergence  & Stability & Robustness   \\
   \midrule
   IGA & $\times$ & $\times$ & $\checkmark$ & $\checkmark$  \\
   GIGA & $d_{\Gamma}$ & $\times$ & $\times$ & $\times$   \\
   SGIGA & $[I-I_b](d_{\Gamma})$ & $\times$ & $\checkmark$ & $\checkmark$   \\
SGIGA-multi & $[I-I_b]\left(d_{\Gamma}\left\{1,s,t\right\}\right)$ & $\times$ &  $\times$ & $\times$ \\
cor.GIGA & $d_{\Gamma}^{\text{mod}}$ & $\times$ &  $\times$ & $\times$\\
   \bottomrule
\end{tabular}
    \caption{A comparison of several different GIGA methods of degree two.}
    \label{differentGIGA1}
\end{table*}

\section{New GIGA enrichment space construction}\label{sec5}

In this section, we provide the detailed construction of SGIGA2 with three main concepts, i.e., modified quasi-interpolation, enriched basis function, projection and orthogonalization method. 
\subsection{{Modified} quasi-interpolation}
{This section develops a modified quasi-interpolation scheme to improve convergence characteristics. A central challenge in quasi-interpolation theory involves ensuring precise polynomial reproduction across the computational domain.} Specifically, if $f$ is a quadratic polynomial, the quasi-interpolation should provide an exact representation. As demonstrated in Section \ref{sec2}, traditional quasi-interpolation methods are effective primarily for smooth functions. Therefore, our primary objective is to develop a quasi-interpolation scheme that efficiently handles functions with interfaces, while minimizing the number of elements that fail to reproduce quadratic polynomials.

For a function $f$ that is $C^0$ along the interface, which can be denoted by {$f=f_0|{\Omega_0^p}+f_1|{\Omega_1^p}$}, where $f_0 \in W^{3, \infty}(\Omega_0^p),f_1 \in W^{3, \infty}(\Omega_1^p)$  and $f_0|_\Gamma = f_1|_\Gamma$. We continuously extend $f_0$ and $f_1$ to the {parameter} domain $\Omega^p$ to obtain functions $\bar{f}_0$ and $\bar{f}_1$ such that
\begin{equation}
  \bar{f}_{r}=f_{r}\quad\text{on} \quad\Omega^p_{r}, \quad \|\bar{f}_r\|_{{W^{3,\infty}(\Omega)}}\leq C\|f_{r}\|_{{W^{3,\infty}(\Omega_{r})}},\quad r=0,1.
\end{equation}

\begin{figure}[htbp]
    \centering
    \subfloat{     
        \includegraphics[width=0.45\linewidth]{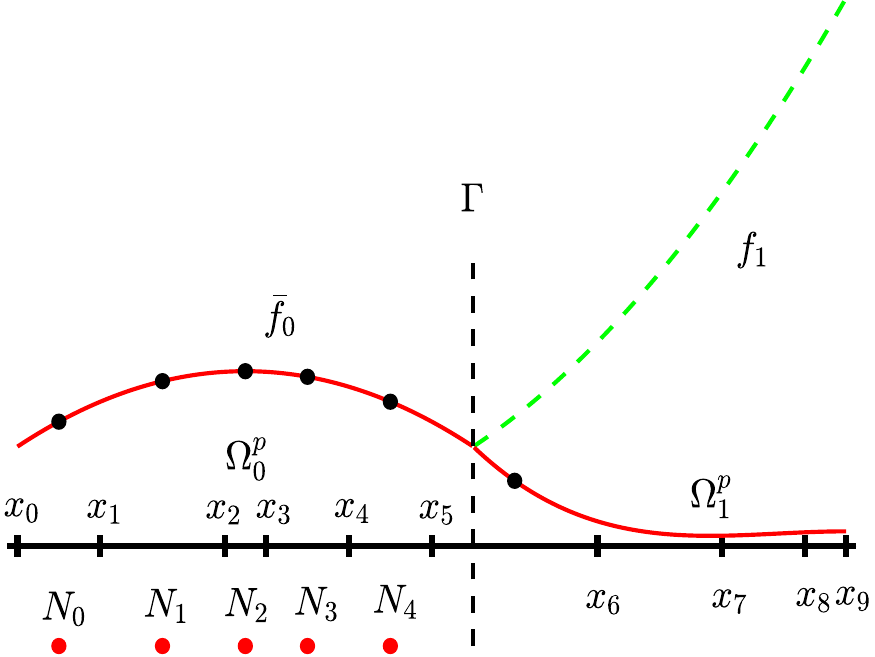}
        }\hfill
    \subfloat{
        \includegraphics[width=0.45\linewidth]{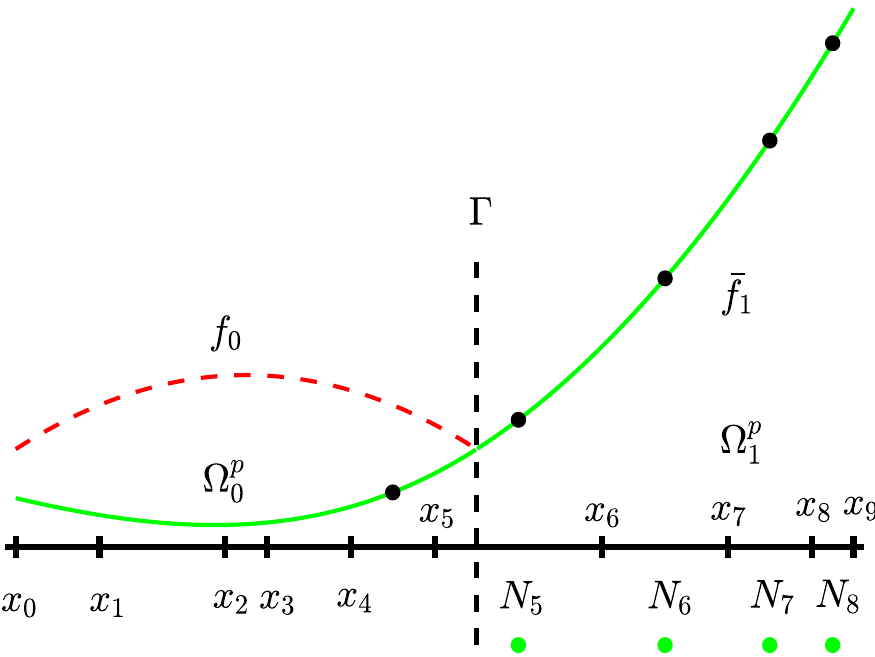}
        }
    \caption{Description of calculation of B-spline coefficients corresponding to modified quasi-interpolation. The solid line shows
the result of the continuation. The red control points indicate that the coefficients corresponding to the B-spline are calculated by the
function represented by the red line. The black dots represent the value used in the calculation.}
    \label{modIb}
\end{figure}

{The basic idea for the modified quasi-interpolation is simple. As illustrated in Fig~\ref{modIb}, for the curve case,
assuming the interface lies within the element $[x_5, x_6]$, then the coefficients of B-splines corresponding to elements in
$\Omega^p_0$ ($N_i,i=0,\dots,4$) are calculated by $\bar{f}_0$ while the coefficients corresponding to elements in
$\Omega^p_1$ and interface element ($N_i,i=5,\dots,8$) are calculated by $\bar{f}_1$.}

We can get the error theorem of the modified quasi-interpolation for the weak discontinuous function.
\begin{theorem}\label{modquasierror}
	There exists a constant $0<C<1$ such that for all {$f=f_0|{\Omega_0^p}+f_1|{\Omega_1^p}$},$f_1 \in W^{3, \infty}(\Omega_1^p),f_2 \in W^{3, \infty}(\Omega_2^p)$  and $f_1|_\Gamma = f_2|_\Gamma$ with step size $h$, we have
	\begin{equation}
	\left\|f-I_b^* f\right\|_{L^{\infty}(\Omega^p/\omega_i)} \leqslant C h^3\left\|f\right\|_{W^{3, \infty}(\Omega^p/\omega_i)}\quad i\in \mathcal{J}^{f}_{1, +}
	\end{equation}
\end{theorem}
The proof of Theorem \ref{modquasierror} follows the same structure as that of Theorem \ref{quasierror}.

The above construction can be extended to 2D case easily,
\begin{equation}
    I_b^*{f} =
    \begin{cases}
    I_b\bar{f}_0, & (s,t)\in \omega_i \subset\Omega^p_0,i\notin \mathcal{J}^{f}_{1, +}, \\ I_b\bar{f}_1, & (s,t)\in \omega_i \subset\Omega^p_1.\\
    \end{cases}
\end{equation}

\begin{figure}[htbp]
    \centering
    \subfloat[The Fig shows the
    result of the traditional quasi-interpolation]{    
        \includegraphics[width=0.4\linewidth]{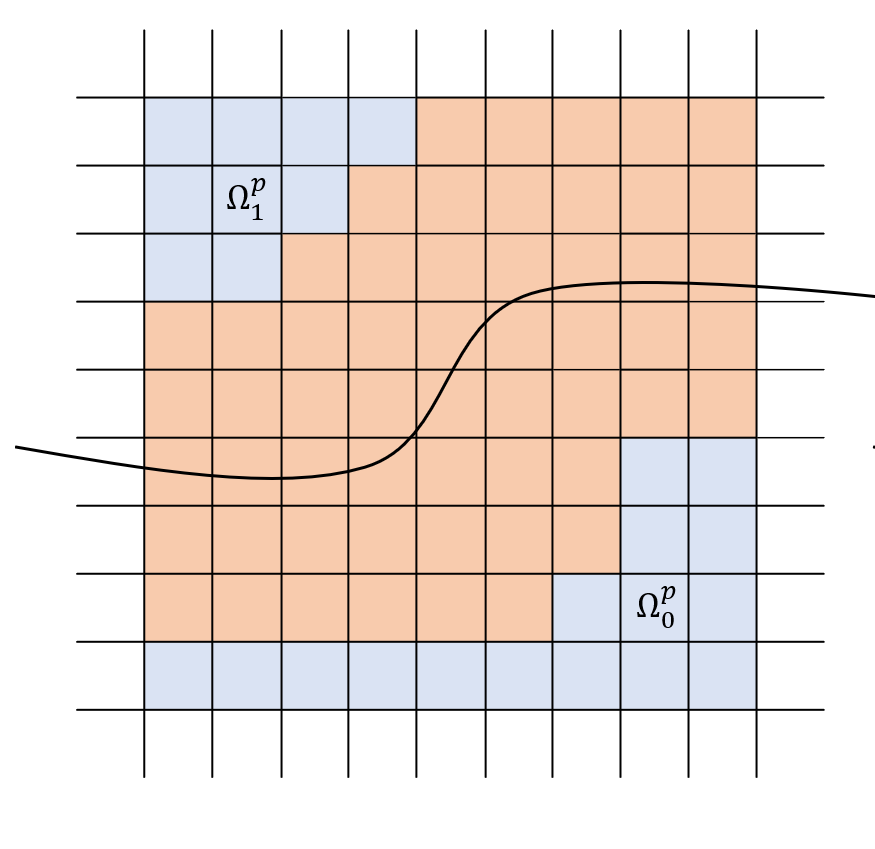}
        }\hfill
    \subfloat[The Fig shows the
    result of the modified quasi-interpolation]{
        \includegraphics[width=0.4\linewidth]{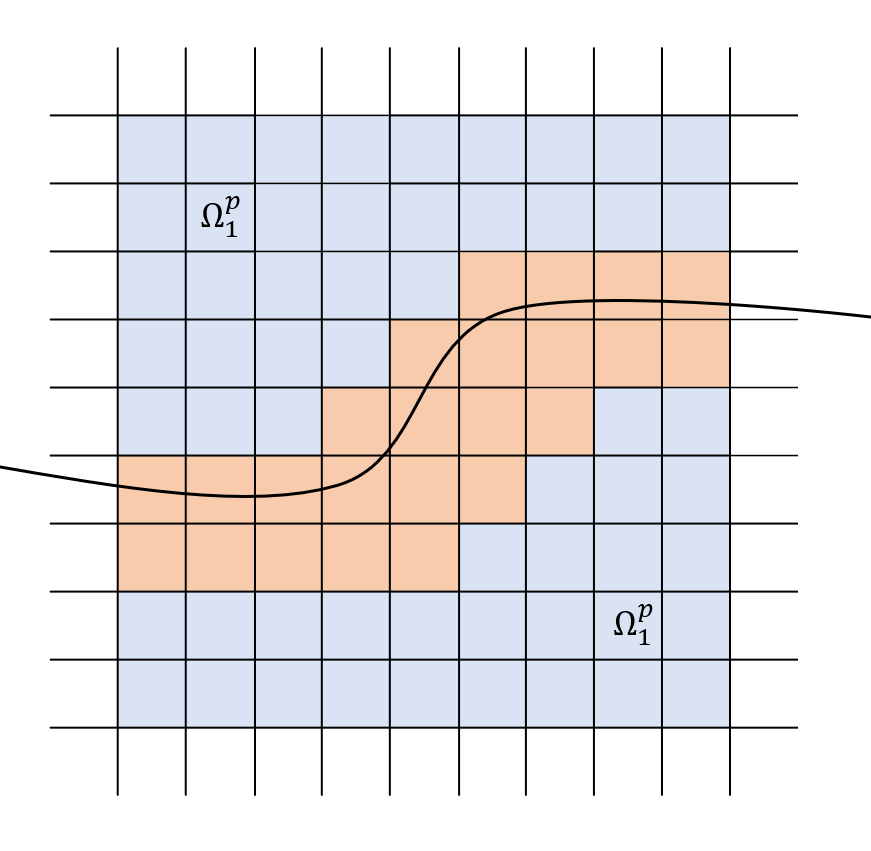}
        }
    \caption{The red elements doe not satisfy the error theorem while the blue elements satisfy the error theorem.}
   \label{fig_2dquasi}
\end{figure}

{Fig~\ref{fig_2dquasi} provides a comparative visualization of approximation accuracy between conventional and enhanced quasi-interpolation schemes. The enhanced formulation extends the applicability of the error estimation theorem to significantly larger element subsets.}

\begin{figure}[htbp]
    \centering
    \subfloat[The logarithm of absolute errors of the initial quasi-interpolation]{     
        \includegraphics[width=0.45\linewidth]{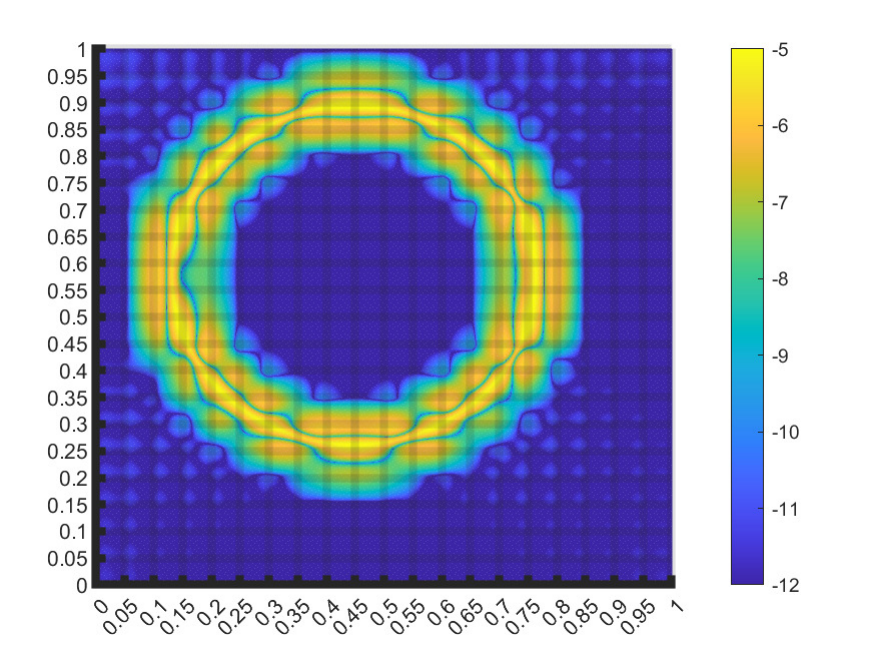}
        }\hfill
    \subfloat[The logarithm of absolute errors of the modified quasi-interpolation]{
        \includegraphics[width=0.45\linewidth]{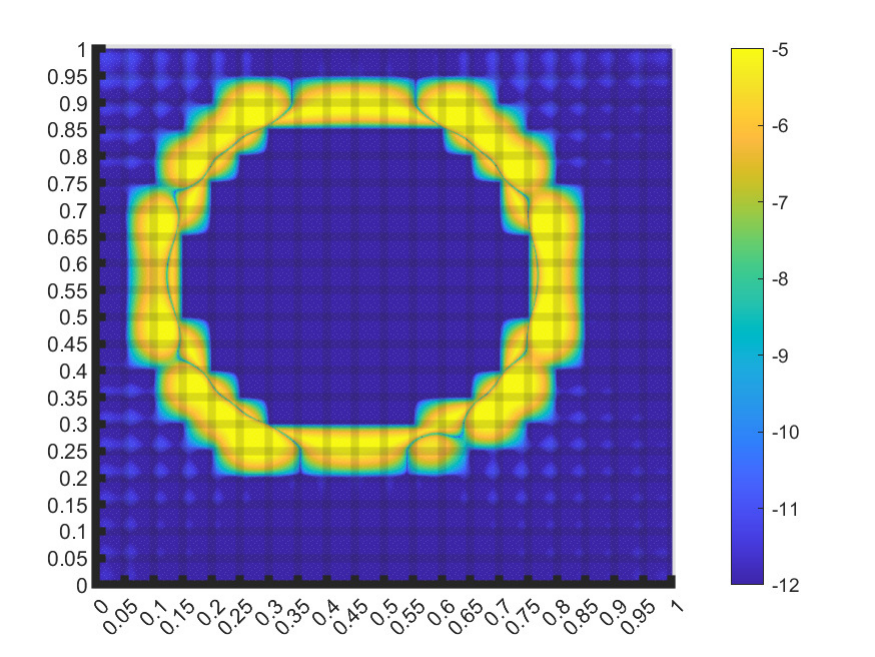}
        }
    \caption{The logarithm of absolute errors of two quasi-interpolations on a function with interface, where blue means the error is zero.}
    \label{mod quasi-interpolation}
\end{figure}

{We also compare the quasi-interpolation for function $f(s, t) = (s-\frac{1}{\sqrt{5}})^2 + (t-\frac{1}{\sqrt{3}})^2- \frac{1}{10}$ if $(s-\frac{1}{\sqrt{5}})^2 + (t-\frac{1}{\sqrt{3}})^2 \geq \frac{1}{10}$ and zeros in the rest region.
Fig \ref{mod quasi-interpolation} shows the results of errors of two quasi-interpolations, where we can {observe} that the modified quasi-interpolation
can reproduce quadratic polynomial in more elements. }

\subsection{New enrichment space and basis function}
\label{sec52}
To improve the convergence rate, we propose SGIGA2, whose enrichment space based on the modified quasi-interpolation and special enriched function $\theta_j$, which is defined on the element of the index set $\mathcal{J}^{f}_{1, +}$. In other words, we define a basis function for each enriched element.

\begin{figure}[htbp]
    \centering
    \subfloat[The Fig illustrates how to compute $\mu_j$]{    
        \includegraphics[width=0.4\linewidth]{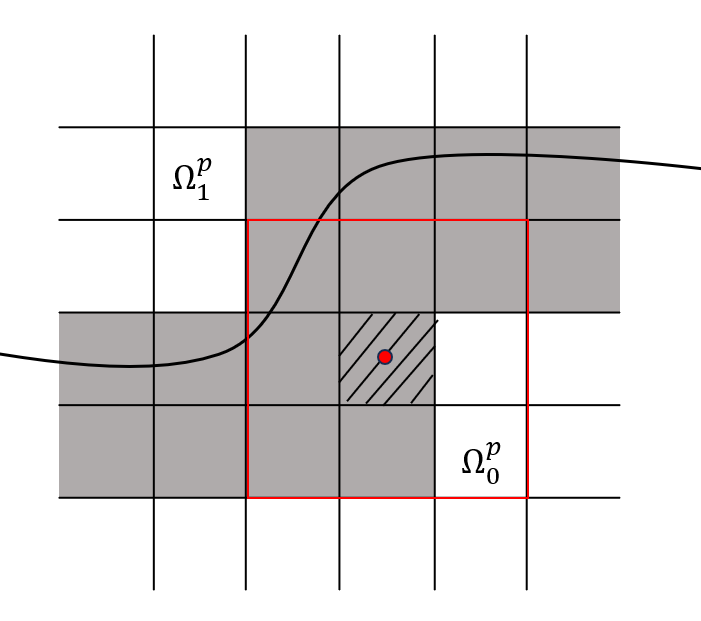}
        }\hfill
    \subfloat[The Fig illustrates how to compute $\theta_i$]{
        \includegraphics[width=0.4\linewidth]{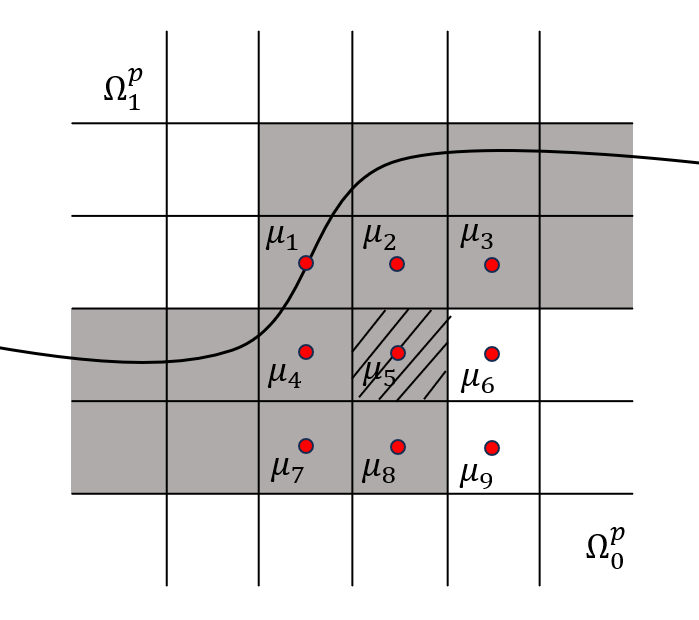}
        }
    \caption{The calculation method of coefficient. And the grey elements are enriched elements.}
   \label{SGIGAenrichment}
\end{figure}

First, for each B-spline basis function $N_{j}(s, t)$, if its support intersects the enriched elements, then denote $\mu_{j}$ to be the number of enriched elements in support of $N_{j}(s, t)$. For example, $\mu_{j}$
for the B-spline basis function corresponding the red control point in Fig~\ref{SGIGAenrichment} (left Fig) is $7$. Now, we are ready to define $\theta_j$ as 
\begin{equation}
\label{PU}
\theta_j = \sum_{j_i\in I_j} \frac{1}{\mu_{j_i}}N_{j_i},
\end{equation}
{where $I_j$ is the index that the support of $N_{j_i}$ intersect the element $j$. It is obvious that $\sum \theta_j$ is one in all the enriched elements, which is related to the proof of the optimal convergence property (equation \eqref{formula36}).}

{We will use Fig~\ref{SGIGAenrichment} to explain the equation. For example, for the shadow element, there are nine
B-spline basis functions contribute the element, each one has a $\mu_{j}$, where $\mu_1={7}$, $\mu_2={8}$, $\mu_3={7}$, $\mu_4={8}$, $\mu_5={7}$, $\mu_6={5}$, $\mu_7={6}$, $\mu_8={4}$, $\mu_9={2}$. These values are used to compute $\theta_i$ for the shadow element.}

To simplify the proof and reduce program overhead, we consider using the unilateral distance function
\begin{equation}
    \tilde{d_\Gamma} =  \begin{cases}d_\Gamma(s,t), & (s,t)\in \Omega^p_0, \\ 0, & (s,t)\in \Omega^p_1 .\end{cases},
\end{equation}
and the approximation space is given by
\begin{equation}
\label{SGIGAmulti}
 \mathbb{V}_{E N R}^h=\operatorname{span}\left\{\theta_j \left(\left(I-I_b^*\right)(\tilde{d_\Gamma}\{1, s, t\})\right): j \in \mathcal{J}^{f}_{1, +} \right\},
\end{equation}
where $I_b^*$ is the modified quasi-interpolation and $\mathcal{J}^{f}_{1, +}$ is the indices set of enriched elements.

We also propose a method to further improve the behavior of the SCN by modifying the basis functions of the enrichment space, which does not change the approximation spaces.
\begin{itemize}
  \item \emph{Local $L^2$ norm projection}: {$T$} is denoted as the $L^2$ norm local orthogonal projection operator. For a given subspace $V_{IGA}^{h*} \subset V_{IGA}^{h}$

  \begin{equation}
  {T}: V_{ENR}^{h} \rightarrow V_{IGA}^{h*}, v \mapsto w=P(v), \end{equation}
  and \begin{equation} (u, v)=\left(u, {T}(v)\right), \forall u \in V_{IGA}^{h*},\end{equation}
  which is equivalent to a linear system $\mathbf{M } {\mathbf{T}}=\mathbf{G}$, where $\mathbf{M}_{i j}=\left(N_j, N_i\right), \mathbf{G}_{i l}=\left(\psi_l, N_i\right)$.
  Thus, {$T$} has a matrix form ${\mathbf{T}}=\mathbf{M}^{-1} \mathbf{G}$,  and the enrichment space is modified to
  \begin{equation}V_{ENR}^{h}\left(T\right)=\operatorname{span}\left\{\psi_l^*:=\psi_l-T\left(\psi_l\right), l \in \mathcal{J}_{h}^E\right\}.\end{equation}
\item \emph{Orthogonalization}: Since $\mathbf{K}_{EE}$ is a positive-definite matrix with a unique decomposition $\mathbf{LDL}^T$, where $\mathbf{D}$ is a diagonal matrix and $\mathbf{L}$ is a full-rank upper triangular matrix. We can define a bijection $L$ on $V_{ENR}^{h}\left(P\right)$. Assume that $\mathbf{X}$ is the coefficient vector of $v$ and $L(v)$ correspond to $\mathbf{L}^{-T} \mathbf{X}$, then we can derived from
    \begin{equation}
    \left(\mathbf{L}^{-T} \mathbf{X}\right)^T \mathbf{K}_{EE}\left(\mathbf{L}^{-T} \mathbf{Y}\right)=\mathbf{X}^T \mathbf{D} \mathbf{Y}
    \end{equation}
    that $L\left(\varphi_l\right), l \in \mathcal{J}_{h}^E$ are orthogonal in the $L^2$ norm.
\end{itemize}

Denoting $I$ to be an identity mapping, the final SGIGA2 enrichment space can be characterized by
\begin{equation}
\begin{aligned}
& V_{E N R}^{h}(S G I G A2):=\operatorname{span}\left\{\psi^*_l:=L\left(I-{T}\right)\left(\psi_l\right), l \in \mathcal{J}_{h}^E\right\} ; \\
& V^{h}=V_{I G A}^{h} \oplus V_{E N R}^{h, p}(S G I G A2) .
\end{aligned}
\end{equation}

\subsection{Proof for the optimal convergence rates} \label{sec51}
{This section derives a error estimates for SGIGA2 approximations of the interface problem. Let $ C>0 $ denote a generic positive constant, dependent on mesh regularity and parametric domain $\Omega^p$ but independent of mesh size $ h $, whose value may vary between occurrences.} Recall that
\begin{equation}
{u}(s,t):= \begin{cases}u_0(s,t), & (s,t) \in \Omega_0^p, \\ u_1(s,t), & (s,t) \in \Omega_1^p .\end{cases}
\end{equation}
Where both $ u_0 $ and $ u_1 $ are smooth. {we construct their continuous extensions $ \tilde{u}_r \in H^3(\Omega^p)\bigcap W^{2,\infty}(\Omega^p) (r = 0,1)$ over the parametric domain $\Omega^p$ satisfying:} 
\begin{equation}
\begin{aligned}
\tilde{u}_r=u_r \text { on } \Omega_r^p \text { and }\left\|\tilde{u}_r\right\|_{H^3(\Omega^p)} \leq C\left\|u_r\right\|_{H^3\left(\Omega_r^p\right)}\\
\quad\left\|\tilde{u}_r\right\|_{W^{2, \infty}(\Omega^p)} \leq C\left\|u_r\right\|_{W^{2, \infty}\left(\Omega_r^p\right)}, \quad r=0,1,
\end{aligned}
\end{equation}
{To facilitate the proof, we define}
\begin{equation}
\tilde{u}(s,t):= \begin{cases}(\tilde{u}_0-\tilde{u}_1)(s,t), & (s,t) \in \Omega^p_0, \\ 0, & (s,t) \in \Omega^p_1 .\end{cases}
\end{equation}
{Two fundamental lemmas essential for establishing the a error estimates of SGIGA2 are formulated as follows:}
\begin{lemma}
For each enriched element $\omega_i$, there is a linear polynomial $\eta_i$ such that
\begin{equation}
\begin{aligned}
 	\left\|\tilde{u}-\tilde{d_\Gamma} \eta_i\right\|_{H^l\left(\omega_i\right)}^2  \leq C h^{6-2 l}\left\|\tilde{u}_0 - \tilde{u}_1 \right\|_{H^3\left(\omega_i\right)}^2+C h^{8-2 l}, l=0,1,
  \end{aligned}
\end{equation}
	\begin{proof}
		According to our previous definition $ \tilde{u}(s,t) = \tilde{d_\Gamma}(s,t) = 0 ,(s,t) \in (\omega_i \bigcap \Omega^p_1)$. Therefore
		\begin{equation}
			\left\|\tilde{u}-\tilde{d_\Gamma} \eta_i\right\|_{H^l\left(\omega_i \bigcap \Omega^p_1 \right)}^2 = 0
		\end{equation}
		In fact, we just need to estimate $ \left\|\tilde{u}-\tilde{d_\Gamma} \eta_i\right\|_{H^l\left(\omega_j \bigcap \Omega^p_0 \right)}^2 $. {Let $ T \tilde{u} $ and $ T \tilde{d_\Gamma} $ denote second-order Taylor polynomials of $ \tilde{u} $ and $ \tilde{d_\Gamma} $ respectively at a point $ (s_\Gamma,t_\Gamma) \in  \omega_i \bigcap \Gamma $ where the local orthogonal coordinate system $ (s,t) $ is equipped with orthonormal basis vectors $ \vec{n}_s $ (tangential to the interface) and $ \vec{n}_t $ (normal to the interface).} And denote $(s,t),(s_\Gamma,t_\Gamma)$ by $p,p_{\Gamma}$. We have

		\begin{equation}	
		\begin{aligned}
		T\tilde{u}(p)&=\tilde{u}(p_{\Gamma}) + {\tilde{u}}_{s}\left(p_{\Gamma}\right)\left(p-p_{\Gamma}\right) \cdot \vec{n}_s +  {\tilde{u}}_{t}\left(p_{\Gamma}\right)\left(p-p_{\Gamma}\right) \cdot \vec{n}_t  +{ \tilde{u}}_{st }\left(p_{\Gamma}\right)\left(\left(p-p_{\Gamma}\right) \cdot \vec{n}_s\left(p-p_{\Gamma}\right) \cdot \vec{n}_t\right)\\
		& +\frac{1}{2} { \tilde{u}}_{tt}\left(p_{\Gamma}\right)\left(\left(p-p_{\Gamma}\right) \cdot \vec{n}_t\right)^2 + \frac{1}{2} { \tilde{u}}_{ ss}\left(p_{\Gamma}\right)\left(\left(p-p_{\Gamma}\right) \cdot \vec{n}_s\right)^2,\quad p \in \omega_i \bigcap \Omega_0^p.
		\end{aligned}
		\end{equation}
        {Observe that $ \tilde{u} $ vanishes identically on $\Gamma$. Consequently, its tangential derivative satisfies $  \tilde{u}_s = \nabla \tilde{u} \cdot \vec{n}_s = 0$,  yielding:}
		\begin{equation}	
		\begin{aligned}
		T\tilde{u}(p)={\tilde{u}}_{t}\left(p_{\Gamma}\right)\left(p-p_{\Gamma}\right) \cdot \vec{n}_t +{ \tilde{u}}_{st }\left(p_{\Gamma}\right)\left(\left(p-p_{\Gamma}\right) \cdot \vec{n}_s\left(p-p_{\Gamma}\right) \cdot \vec{n}_t\right)+
		+
		\frac{1}{2} { \tilde{u}}_{ tt}\left(p_{\Gamma}\right)\left(\left(p-p_{\Gamma}\right) \cdot \vec{n}_t\right)^2,
		\end{aligned}
		\end{equation}
		Similarly,
		\begin{equation}
		\begin{aligned}
		T d_\Gamma(p)&=  { d_\Gamma}_{s}\left(p_{\Gamma}\right)\left(p-p_{\Gamma}\right) \cdot \vec{n}_s+{ d_\Gamma}_{t}\left(p_{\Gamma}\right)\left(p-p_{\Gamma}\right) \cdot \vec{n}_t+\frac{1}{2} { d_\Gamma}_{ss}\left(p_{\Gamma}\right)\left(\left(p-p_{\Gamma}\right) \cdot \vec{n}_s\right)^2 \\
		& +{ d_\Gamma}_{st}\left(p_{\Gamma}\right)\left(\left(p-p_{\Gamma}\right) \cdot \vec{n}_{s}\left(p-p_{\Gamma}\right) \cdot \vec{n}_{t}\right)+\frac{1}{2} { d_\Gamma}_{tt}\left(p_{\Gamma}\right)\left(\left(p-p_{\Gamma}\right) \cdot \vec{n}_{t}\right)^2 \\
		&=  \left(p-p_{\Gamma}\right) \cdot \vec{n}_t\left(1+\frac{1}{2} { \tilde{d_\Gamma}}_{tt}\left(p_{\Gamma}\right)\left(p-p_{\Gamma}\right) \cdot \vec{n}_t\right).
		\end{aligned}
		\end{equation}
		Therefore
		\begin{equation}
		\begin{aligned}
		&\left\|\tilde{u}-\tilde{d_\Gamma} \eta_i\right\|_{H^l\left(\omega_i \bigcap \Omega^p_0 \right)}^2 \\&=
		\left\|\tilde{u} - T\tilde{u} + T\tilde{u} - T\tilde{d_\Gamma}\eta_i - (\tilde{d_\Gamma} - T\tilde{d_\Gamma})\eta_i \right\|_{H^l\left(\omega_i \cap \Omega^p_0\right)}^2\\
		& \leq \left\|\tilde{u} - T\tilde{u}\right\|_{H^l\left(\omega_i \bigcap \Omega^p_0 \right)}^2 + \left\|T\tilde{u} - T\tilde{d_\Gamma}\eta_i\right\|_{H^l\left(\omega_i \bigcap \Omega^p_0 \right)}^2
		+ \left\|(\tilde{d_\Gamma} - T\tilde{d_\Gamma})\eta_i\right\|_{H^l\left(\omega_i \bigcap \Omega^p_0 \right)}^2\\
		&=: A_1 + A_2 + A_3.
		\end{aligned}
		\end{equation}
	where $ A_1 $ and $ A_3 $ can be estimated by Taylor expansion:
	\begin{equation}\label{a1}
	\begin{aligned}
	A_1 = &\left\|\tilde{u} - T\tilde{u}\right\|_{H^l\left(\omega_i \bigcap \Omega^p_0 \right)}^2 \leqslant \left\|\tilde{u} - T\tilde{u}\right\|_{H^l\left(\omega_i \right)}^2
	\leqslant Ch^{6-2l} \left\|\tilde{u}_0-\tilde{u}_1\right\| _{H^l\left(\omega_i \right)}^2\\
	 A_3 = &\left\|(\tilde{d_\Gamma} - T\tilde{d_\Gamma})\eta_i\right\|_{H^l\left(\omega_i \bigcap \Omega^p_0 \right)}^2 \leqslant \left\|(\tilde{d_\Gamma} - T\tilde{d_\Gamma})\right\|_{W^{l,\infty}\left(\omega_i\right)}^2 \left\|\eta_i\right\|_{W^{l,\infty}\left(\omega_i\right)}^2 |\omega_i|^2  \leqslant Ch^{8-2l}.
	 \end{aligned}
	\end{equation}
    To estimate $ A_2 $, let $ T\tilde{u} = (p-p_{\Gamma})\cdot \vec{n}_t (c_0 + c_1 (p-p_{\Gamma})\cdot \vec{n}_t + c_2 (p-p_{\Gamma})\cdot \vec{n}_s)$ and $ \eta_i = d_0 + d_1(p-p_{\Gamma})\cdot \vec{n}_t+ d_2 (p-p_{\Gamma})\cdot \vec{n}_s $, where $ c_j,j=0,1,2$
	are determined by the Taylor expansion of $ \tilde{u} $ while $ d_j,j=0,1,2$ are undetermined.
\begin{small}
	\begin{equation}
	\begin{aligned}
	A_2 &= \left\|T\tilde{u} - T\tilde{d_\Gamma}\eta_i\right\|_{H^l\left(\omega_i \bigcap \Omega^p_0 \right)}^2 \\
	&= \left\| \left(p-p_{\Gamma}\right)  \vec{n}_t\left(c_0 + c_1 (p-p_{\Gamma})\cdot \vec{n}_t + c_2 (p-p_{\Gamma}) \vec{n}_s - (1+\frac{1}{2} { \tilde{d_\Gamma}}_{tt}\left(p_{\Gamma}\right)\left(p-p_{\Gamma}\right)  \vec{n}_t) \eta_i \right) \right\|_{H^l\left(\omega_i \bigcap \Omega^p_0 \right)}^2.
	\end{aligned}		
	\end{equation}
\end{small}
	Then we can determine $ d_0 = c_0,d_1 = c_1 - d_0 \frac{1}{2} { \tilde{d_\Gamma}}_{tt},d_2 = c_2 $ and denote $ \frac{1}{2} { \tilde{d_\Gamma}}_{tt} $ by $ d $.

	\begin{equation}\label{a2}
	\begin{aligned}
	A_2
	&= \left\| \left(p-p_{\Gamma}\right) \cdot \vec{n}_t\left( d_1 d (\left(p-p_{\Gamma}\right) \cdot \vec{n}_t)^2 + d_2 d (\left(p-p_{\Gamma}\right) \cdot \vec{n}_t)\cdot (\left(p-p_{\Gamma}\right) \cdot \vec{n}_s) \right) \right\|_{H^l\left(\omega_i \bigcap \Omega^p_0 \right)}^2\\
	&\leq \left\| \left( d_1 d (\left(p-p_{\Gamma}\right) \cdot \vec{n}_t)^3 + d_2 d (\left(p-p_{\Gamma}\right) \cdot \vec{n}_t)^2\cdot (\left(p-p_{\Gamma}\right) \cdot \vec{n}_s) \right) \right\|_{W^{l,\infty}\left(\omega_i \bigcap \Omega^p_0 \right)}^2 \cdot |\omega_i|^2\\
	& \leq Ch^{8-2l}.
	\end{aligned}		
	\end{equation}
	{Combining \eqref{a1} with \eqref{a2} yields:}
	\begin{equation}
		\begin{aligned}
		\left\|\tilde{u}-\tilde{d_\Gamma} \eta_i\right\|_{H^l\left(\omega_i\right)}^2 &= \left\|\tilde{u}-\tilde{d_\Gamma} \eta_i\right\|_{H^l\left(\omega_i \cap \Omega^p_0 \right)}^2 \leq A_1 + A_2 + A_3\\
		&\leq C h^{6-2 l}\left\|\tilde{u}_0 - \tilde{u}_1 \right\|_{H^3\left(\omega_i\right)}^2+C h^{8- 2l}, l=0,1.
		\end{aligned}
	\end{equation}

	\end{proof}

\end{lemma}

\begin{lemma}
	For each enriched element $\omega_i$, there is a linear polynomial $\eta_i$ such that
	\begin{equation}
     \left\|I_b^*(\tilde{u}-\tilde{d_\Gamma} \eta_i)\right\|_{H^l\left(\omega_i\right)}^2  \leq C h^{6-2 l}\left\|\tilde{u}_0 - \tilde{u}_1 \right\|_{H^3\left(\omega_i\right)}^2+C h^{8-2 l}, l=0,1,
	\end{equation}
	\begin{proof}
		where $ I_b^* $ is modified quasi-interpolation. We have
		\begin{equation}
		\begin{aligned}
		I_b^*(\tilde{u}-\tilde{d_\Gamma} \eta_i) &= \sum_{k,j} (\tilde{u}-\tilde{d_\Gamma} \eta_i)(\tau_j^k) N_j \\
		&= \sum_{k,j} (\tilde{u}-T\tilde{u} )(\tau_j^k) N_k + \sum_{k,j} (T\tilde{u}-T\tilde{d_\Gamma}\eta_i)(\tau_j^k) N_k
  + \sum_{k,j} ((\tilde{d_\Gamma}-T\tilde{d_\Gamma})\eta_i)(\tau_j^k) N_k
		\end{aligned}
		\end{equation}
		
		{Applying Taylor's theorem with remainder to \eqref{a1} and \eqref{a2} establishes}
		\begin{equation}
		\begin{aligned}
		\left\|{I}_b^*(\tilde{u}-\tilde{D} \eta_i)\right\|_{H^l\left(\omega_i\right)}^2 &\leq   \left\{ C\max_{\substack{\tau_j^k }}\left|\left(\tilde{u}-T \tilde{u}\right)\left(\tau_j^k\right)\right|^2 + Ch^6\right\} \sum_k\left\|N_k\right\|_{H^l\left(\omega_i\right)}^2  \\
		&\leq  Ch^{6-2l}\left\|\tilde{u}_0 - \tilde{u}_1 \right\|_{H^3\left(\omega_i\right)}^2+C h^{8-2 l}
		\end{aligned}
		\end{equation}
	\end{proof}
\end{lemma}

\begin{theorem}\label{erroranalysis}
	(Error Estimate for SGIGA2). Let $u\in C^0(\Omega^p)$ be the solution of problem \eqref{problem} with $ u|\Omega^p_i \in H^3(\Omega^p_i)$
	for $ i = 0, 1 $, and {let $u_h$ denote the SGIGA2 approximation.} Then
	\begin{equation}	
	|u-u_h|_{H^1\left(\Omega^p\right)} \leq Ch^2 \|u\|_{H^3\left(\Omega^p\right)}.
	\end{equation}
	\begin{proof}
    {The proof construction requires demonstrating the existence of $ u_h^* \in \mathbb{V}_{h} $ fulfilling the approximation criterion:}
		\begin{equation}
		|u-u_h|_{H^1\left(\Omega^p\right)} \leq |u-u_h^*|_{H^1\left(\Omega^p\right)}\leq Ch^2 \|u\|_{H^3\left(\Omega^p\right)}.
		\end{equation}
Define $ u_h^* = I_b^*\tilde{u}+I_b^*\tilde{u_1}+ v$, it is clear that $ I_b^*\tilde{u},I_b^*\tilde{u}_1\in \mathbb{V}_{I G A}^h $ and $ v \in  \mathbb{V}_{ENR}^h$. $ v $ can be expressed as
\begin{equation}
v:=\sum_{j\in\mathcal{J}_h^{E}}\theta_j\left(\tilde{d_\Gamma}\eta_j - {I}_b^*(\tilde{d_\Gamma}\eta_j)\right),
\end{equation}
{where $ \eta_j $ denotes a linear polynomial by coefficients associated with the $ \{1, s, t\} $.} Combining these results establishes the theorem through:
		\begin{equation}
		\begin{aligned}
		|u-u_h^*|^2_{H^1\left(\Omega^p\right)} &= |\tilde{u}-I_b^*\tilde{u}+ \tilde{u}_1-I_b^*\tilde{u} _1-v|^2_{H^1\left(\Omega^p\right)}\\
		&\leq |\tilde{u}-I_b^*\tilde{u} - v|^2_{H^1\left(\Omega^p\right)} + |\tilde{u}_1-I_b^*\tilde{u}_1 |^2_{H^1\left(\Omega^p\right)}
		= :T_1 + T_2.
		\end{aligned}
		\end{equation}
		 The error estimate for the $ T_2 $  can be obtained directly from the error estimate of the quasi-interpolation 		
   \begin{equation}
		 T_2 = |\tilde{u}_1-I_b^*\tilde{u}_1 |^2_{H^1\left(\Omega^p\right)} \leq Ch^4\|\tilde{u}_1\|^2_{H^3\left(\Omega^p\right)} \leq Ch^4 \|{u}\|^2_{H^3\left(\Omega^p\right)}.
  \end{equation}
		 And $ T_1 = |\tilde{u}-I_b^*\tilde{u} - v|^2_{H^1\left(\Omega^p\right)}  = \sum |\tilde{u}-I_b^*\tilde{u} - v|^2_{H^1\left(\omega_i\right)}$. For the estimation of $ T_1 $, all elements can be divided into three categories: (1) enrichment element $ \mathcal{T}_e = \mathcal{J}^{f}_{1, +}$, (2) blending element $\mathcal{T}_b = \mathcal{J}^{f}_{3, +} \setminus \mathcal{J}^{f}_{1, +} \cup \mathcal{J}^{f}_{2, -} \setminus \mathcal{J}^{f}$, and (3) ordinary element $\mathcal{T}_n$. Thus,
		\begin{equation}
		\begin{aligned}
		T_1 &= \sum_{i \in \mathcal{T}_e} |\tilde{u}-I_b^*\tilde{u} - v|^2_{H^1\left(\omega_i\right)} + \sum_{i \in \mathcal{T}_b} |\tilde{u}-I_b^*\tilde{u} - v|^2_{H^1\left(\omega_i\right)} + \sum_{i \in \mathcal{T}_n} |\tilde{u}-I_b^*\tilde{u} |^2_{H^1\left(\omega_i\right)}\\
		& =:A_1 + A_2 + A_3
		\end{aligned}
		\end{equation}
        {Let $ \mathcal{J}_{h,i}^{E} $ denote the enriched basis function index set for element $ \omega_i $. The estimation of $A_3$ follows directly via the quasi-interpolation approximation property:}
	\begin{equation}
	A_3 = \sum_{i \in \mathcal{T}_n} |\tilde{u}-I_b^*\tilde{u} |^2_{H^1\left(\omega_i\right)}\leq C h^4 \sum_{i \in \mathcal{T}_n} \|\tilde{u}\|_{H^3\left(\omega_i\right)}^2 \leq C h^4\|\tilde{u}\|_{H^3\left(\Omega^p\right)}^2.
	\end{equation}
{For the error estimation of $A_{1}$, we have the following results}
\begin{equation}
		\begin{aligned}
  \label{formula36}
		A_1 &= \sum_{i \in \mathcal{T}_e} |\tilde{u}-I_b^*\tilde{u} - v|^2_{H^1\left(\omega_i\right)} = \sum_{i \in \mathcal{T}_e} |\tilde{u}-I_b^*\tilde{u} - \sum_{j\in\mathcal{J}_{h,i}^{E}}\theta_j\left(\tilde{d_\Gamma}\eta_j - {I}_b^*(\tilde{d_\Gamma}\eta_j)\right)|^2_{H^1\left(\omega_i\right)}\\
		& = \sum_{i \in \mathcal{T}_e} |\sum_{j\in\mathcal{J}_{h,i}^{E}}\theta_j(\tilde{u}-I_b^*\tilde{U}) - \sum_{j\in\mathcal{J}_{h,i}^{E}}\theta_j\left(\tilde{d_\Gamma}\eta_j - {I}_b^*(\tilde{d_\Gamma}\eta_j)\right)|^2_{H^1\left(\omega_i\right)}\\
		& = \sum_{i \in \mathcal{T}_e} |\sum_{j\in\mathcal{J}_{h,i}^{E}}\theta_j(\tilde{u}-\tilde{d_\Gamma}\eta_j) - \theta_j I_b^*(\tilde{u}-\tilde{d_\Gamma}\eta_j)|^2_{H^1\left(\omega_i\right)}\\
		& \leq C \sum_{i \in \mathcal{T}_e}\sum_{j\in\mathcal{J}_{h,i}^{E}}\left[\left|\left(\tilde{u}-\tilde{d_\Gamma}\eta_j\right)\right|_{H^1\left(\omega_i\right)}^2+\left|I_b^*\left(\tilde{u}-\tilde{d_\Gamma}\eta_j\right)\right|_{H^1\left(\omega_i\right)}^2\right] \\
		& \leq C \sum_{i \in \mathcal{T}_e}\left[ h^{4}\left\|\tilde{u}_0 - \tilde{u}_1 \right\|_{H^3\left(\omega_i\right)}^2+ h^{6}
		\right] \\
		& \leq Ch^4\left\|\tilde{u}_0 - \tilde{u}_1 \right\|_{H^3\left(\Omega^p\right)}^2+ Ch^{5}  \leq Ch^4 \|u\|_{H^3\left(\Omega^p\right)},
		\end{aligned}
		\end{equation}
where the last two inequalities are derived from the previous lemma.
{The estimation of $A_2$ follows through analogous arguments:}
		\begin{equation}
		\begin{aligned}
  \label{formula37}
		A_2 &= \sum_{i \in \mathcal{T}_b} |\tilde{u}-I_b^*\tilde{u} - v|^2_{H^1\left(\omega_i\right)} = \sum_{i \in \mathcal{T}_b} |\tilde{u}-I_b^*\tilde{u} - \sum_{j\in\mathcal{J}_{h,i}^{E}}\theta_j\left(\tilde{d_\Gamma}\eta_j - {I}_b^*(\tilde{d_\Gamma}\eta_j)\right)|^2_{H^1\left(\omega_i\right)}\\
		&\leq \sum_{i \in \mathcal{T}_b} |\tilde{u}-I_b^*\tilde{u}|^2_{H^1\left(\omega_i\right)} + \\
  &C \sum_{i \in \mathcal{T}_b}\sum_{j\in\mathcal{J}_{h,i}^{E}} \left[\left\|\theta_j\right\|_{L^{\infty}(\Omega^p)}^2\left|\tilde{d_\Gamma}\eta_j-{I}_b^* \tilde{d_\Gamma}\eta_j\right|_{H^1\left(\omega_i\right)}^2+\left\|\nabla \theta_j\right\|_{L^{\infty}(\Omega)}^2\left\|\tilde{d_\Gamma}\eta_j-{I}_b^* \tilde{d_\Gamma}\eta_j\right\|_{L^2\left(\omega_i\right)}^2\right]\\
		&\leq C \left[ h^{4}\sum_{i \in \mathcal{T}_b} \|\tilde{u}\|_{H^3\left(\omega_i\right)}^2 + \sum_{i \in \mathcal{T}_b} \left(\left|\tilde{d_\Gamma}\eta_i-{I}_b^* \tilde{d_\Gamma}\eta_i\right|_{H^1\left(\omega_i\right)}^2 + h^{-2} \left\|\tilde{d_\Gamma}\eta_i-{I}_b^* \tilde{d_\Gamma}\eta_i\right\|_{L^2\left(\omega_i\right)}^2\right)
		\right]	\\
		&\leq C \left[h^{4}\|\tilde{u}\|_{H^3\left(\Omega\right)}^2 + h^4 \sum_{i \in \mathcal{T}_b} \left|\tilde{d_\Gamma}\eta_i-{I}_b^*\tilde{d_\Gamma}\eta_i\right|_{H^3\left(\omega_i\right)}^2 \right],
		\end{aligned}
		\end{equation}
		where the $ \sum_{i \in \mathcal{T}_b} \left|\tilde{d_\Gamma}\eta_i-{I}_b^* \tilde{d_\Gamma}\eta_i\right|_{H^3\left(\omega_i\right)}^2 $ is bounded. For the sake of simplicity, we denote $ A_2 \leq C h^{4}\|\tilde{u}\|_{H^3\left(\Omega^p\right)}^2$. According to a few error estimates, we have $T_1 \leq Ch^4 \|\tilde{u}\|_{H^3\left(\Omega^p\right)}^2$.	Finally,
		\begin{equation}
		|u-u_h|_{H^1\left(\Omega^p\right)} \leq |u-u_h^*|_{H^1\left(\Omega^p\right)}\leq Ch^2 \|u\|_{H^3\left(\Omega^p\right)}
		\end{equation}
		
	\end{proof}
\end{theorem}

\section{Numerical experiment}\label{sec6}

{This section validates the theoretical framework through two-dimensional benchmark problems featuring analytical solutions.}
The exact manufactured solution $ u $ of the problem \eqref{problem} will be used, while the flux $ q $ and $ g $ can be calculated.
We discretize the domain $ \Omega^p $ into $ N\times N $  square uniform elements with mesh parameter $ h = \frac{1}{N} $. {Several representative interface configurations are considered to assess: (a) scaled condition number (SCN) growth rates, and (b) $ H^1 $ seminorm convergence behavior across methodologies:}

$ \blacksquare $ IGA, typical IGA space of degree two.

$ \blacksquare $ GIGA, employed the distance function as enrichment. The approximation space is defined in Section \ref{sec4}. The enriched function is B-spline with indices $\mathcal{J}^{f}_{1}$.

$ \blacksquare $ GIGA*, employed $ \tilde{d_\Gamma} - I_b^* \tilde{d_\Gamma} $ as enrichment. The enriched function $ \theta_j $ is constructed by \eqref{PU} with indices $\mathcal{J}^{f}_{1,+}$.

$ \blacksquare $ SGIGA2, employed $ \tilde{d_\Gamma}\left\{1,x,y\right\} - I_b^* \tilde{d_\Gamma} \left\{1,x,y\right\} $ as enrichment. The approximation space is defined in \eqref{SGIGAmulti}, the enriched function is $ \theta_j $ with indices $\mathcal{J}^{f}_{1,+}$. {To maintain numerical robustness and stabilize the solution, projection and orthogonalization are implemented, enforcing SCN growth rates compatible with standard IGA space.}

\subsection{{An example of a circular interface}}
Let the computational domain be defined as $ \Omega := [0,1] \times [0,1] $ , which is divided into two parts by the circular interface. The circular interface has an equation $ (x-\frac{1}{\sqrt{5}})^2 + (y-\frac{1}{\sqrt{3}})^2= \frac{1}{10} $. Mark the area inside the interface as
$ \Omega_0 $, and the other area is $ \Omega_1 $. {The manufactured solution with prescribed interfacial continuity is constructed as}
\begin{equation}
\begin{aligned}
	 \begin{cases} u_0 = \frac{2a_1}{(a_1-a_0)r_0^4}r^2\cos(2\theta), \text { if } (x,y)\in \Omega_0 \\
	u_1 = \frac{a_1+a_0}{(a_1-a_0)r_0^4}r^2\cos(2\theta) + r^{-2} \cos(2\theta),\\ \text { if } (x,y)\in \Omega_1
	\end{cases}
\end{aligned}
\end{equation}
where the centre of polar coordinates $ (r,\theta) $ is $ (\frac{1}{\sqrt{5}},\frac{1}{\sqrt{3}}) $ and $ r_0 = \frac{1}{\sqrt{10}} $. This solution is $ C^0 $ but not $ C^1 $ at the interface and the solution is smooth on $ \Omega_0 $, $ \Omega_1 $. {The domain is uniformly discretized into $N \times N$ elements with $ N = 5, 10, 20, \dots , 160 $. This configuration ensures the interface remains embedded within element interiors, constituting a non-conforming discretization where no mesh alignment with $\Gamma$ is required.}

\begin{figure}[htbp]
    \centering
    \subfloat[The $ H^1 $ seminorm errors]{     
        \includegraphics[width=0.45\linewidth]{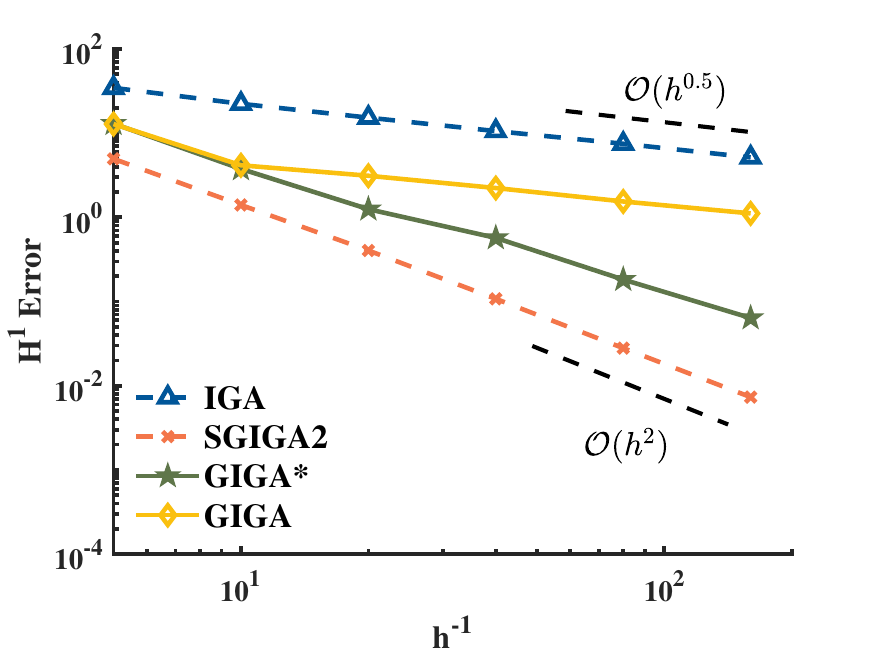}
        }\hfill
    \subfloat[The SCNs results]{
        \includegraphics[width=0.45\linewidth]{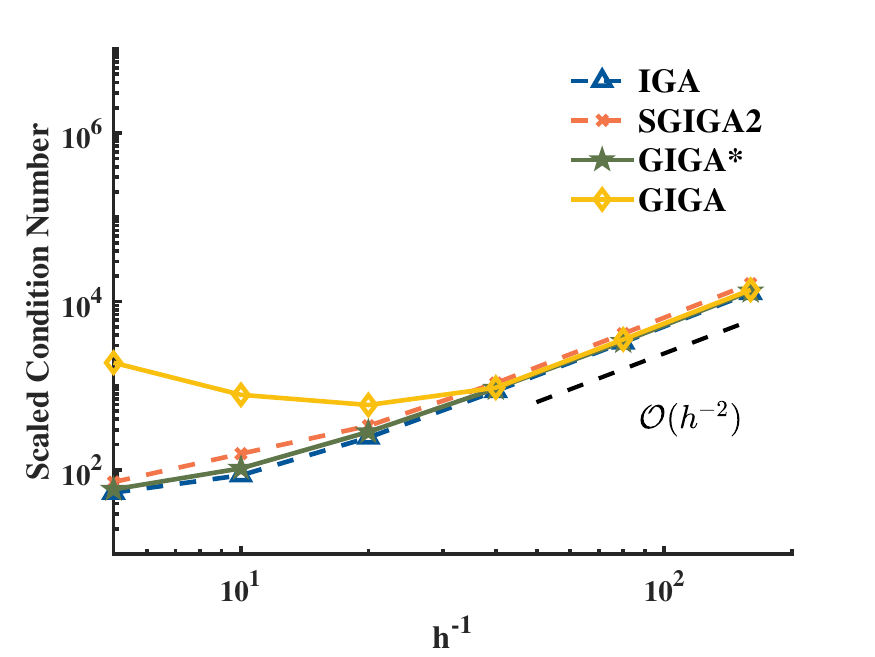}
        }
    \caption{The $ H^1 $ seminorm errors and SCNs with respect to $ h $(from $ \frac{1}{5} $ to $ \frac{1}{160} $) for circular interface($ a_0 = 10, a_1 = 1$).}
    \label{the error of circular interface}
\end{figure}

Fig \ref{the error of circular interface} shows the convergence orders of several methods. The SGIGA2 yields optimal convergence while  the $ H^1 $ seminorm
errors are $ O(h^2) $ as discussed in Theorem \ref{erroranalysis}. {In this example, the convergence rate of GIGA* is much faster than that of the straight line interface example mentioned above, which employs the modified quasi-interpolation $I_b^*$ and the new enriched element selection strategy.} The GIGA method still only reduces the error, but can not improve the rate of convergence. The SCNs results are shown in Fig~\ref{the error of circular interface}. The SCNs of all methods increase with the same order $ O(h^{-2}) $, which is the typical order of the SCN of IGA.

To further study our SGIGA2 method, we demonstrate its application to the ring domain with an arc interface.
\subsection{{An example of an arc interface}}
\label{sec71}
Let the computational domain be defined as $ \Omega  := 1<x^2 + y^2<2$ , which is divided into two parts by the arc interface.  The corresponding parameter $(s,t)$ area $ \Omega^p $ is $[1,2] \times [0,2\pi]$. We consider the mapping relationship where the arc interface corresponds to a straight line in the parameter domain.

\begin{figure}[htbp]
    \centering
    \subfloat[$u$ in the parameter domain]{     
        \includegraphics[width=0.45\linewidth]{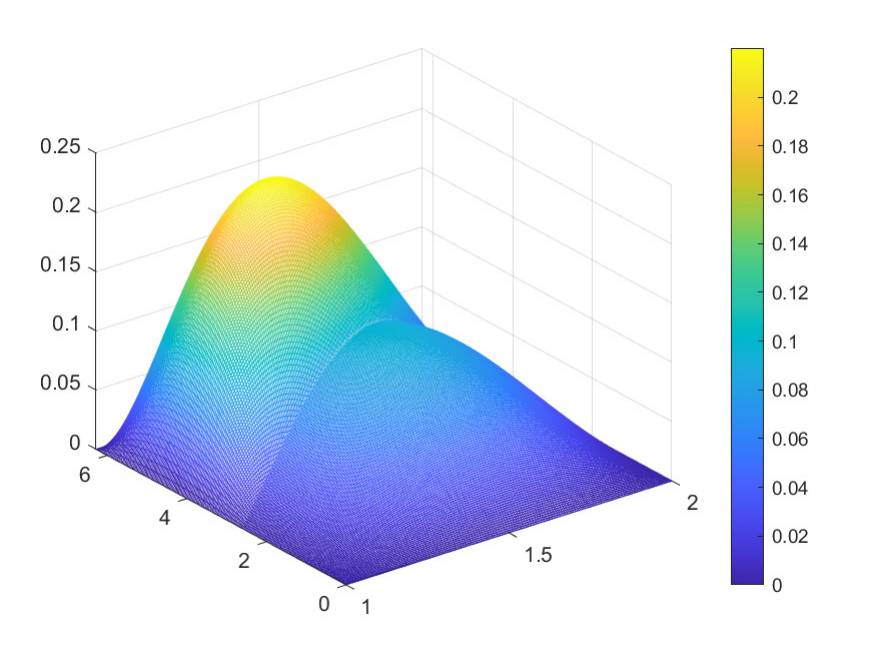}
        }\hfill
    \subfloat[$u$ in the physical domain]{
        \includegraphics[width=0.45\linewidth]{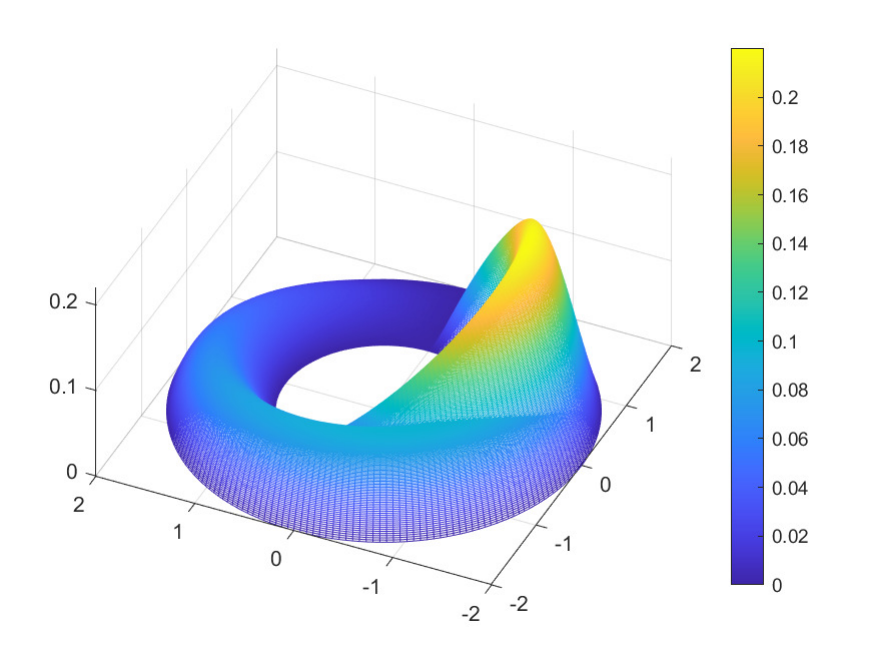}
        }
    \caption{The figure of manufactured solution($a_0 = 20, a_1 = 1$).}
\end{figure}

The straight line in the parameter domain has an equation $ t = 2\pi\tan(\frac{\pi}{8})s $. The manufactured solution of \eqref{problem} in the parameter domain is as follows:
\begin{equation}
   \begin{aligned}
\begin{cases} u_0 = (s-2)(s-1)t(t-2\pi)(r\cos(\theta  - \theta_0) \\+ \frac{a_0}{a_1}r\sin(\theta - \theta_0 )), \text { if } t > 2\pi\tan(\frac{\pi}{8})s (\Omega^p_0) \\
u_1 = (s-2)(s-1)t(t-2\pi)(r\cos(\theta  - \theta_0) \\+ r\sin(\theta  - \theta_0)),  \text { if } t < 2\pi\tan(\frac{\pi}{8})s (\Omega^p_1)
\end{cases}
\end{aligned} 
\end{equation}

where the centre of polar coordinates $ (r,\theta) $ is $ (0,0) $ and $\theta_0 = \arctan(2\pi \tan{\frac{\pi}{8}})$. It's easy to verify that this solution $ u $ is smooth on $ \Omega^p_0 $, $ \Omega^p_1 $ and $ C^0 $ but not $ C^1 $ at the interface.  {The domain is uniformly discretized into $N \times N$ elements with $ N = 5, 10, 20, \dots , 160 $}.

\begin{figure}[htbp]
    \centering
    \subfloat[$ a_0 = 20, a_1 = 1$]{     
        \includegraphics[width=0.45\linewidth]{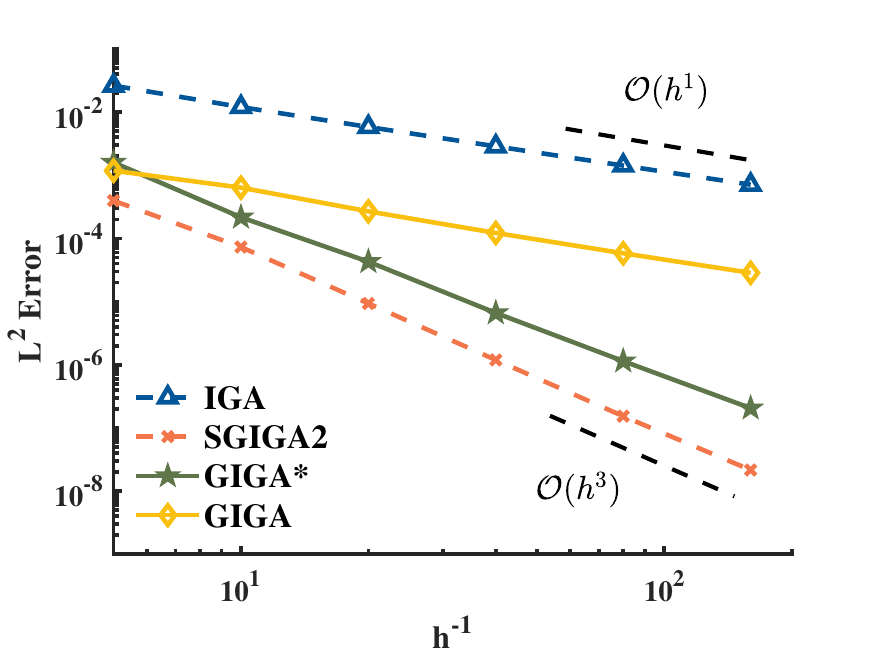}
        }\hfill
    \subfloat[$ a_0 = 100, a_1 = 1$]{
        \includegraphics[width=0.45\linewidth]{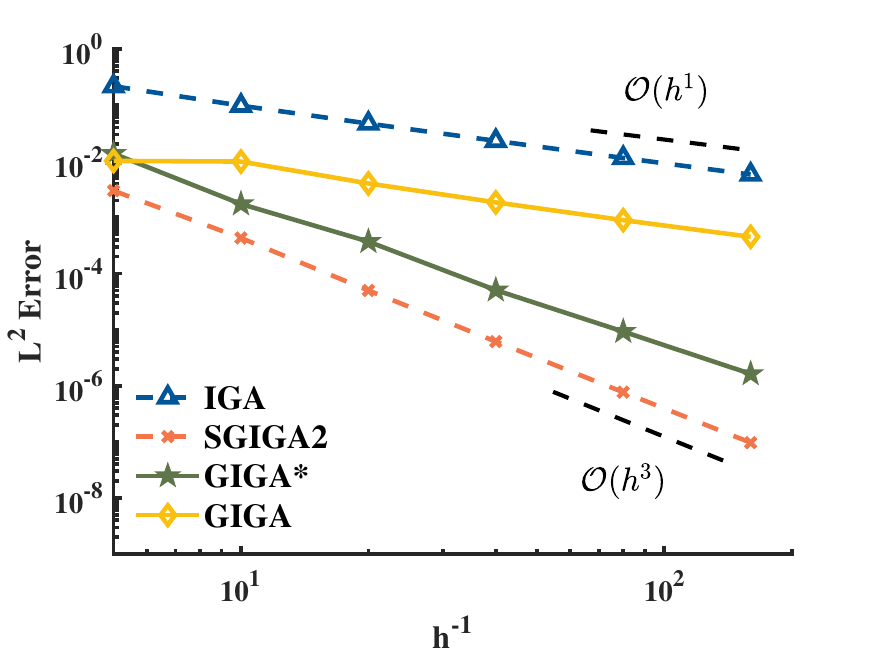}
        }
    \caption{$L^2$ errors with respect to $h$ (from $ \frac{1}{5} $ to $\frac{1}{160}$) for an arc interface.}
    \label{the L2 error of yuanpan interface}
\end{figure}

\begin{figure}[htbp]
    \centering
    \subfloat[$ a_0 = 20, a_1 = 1$]{     
        \includegraphics[width=0.45\linewidth]{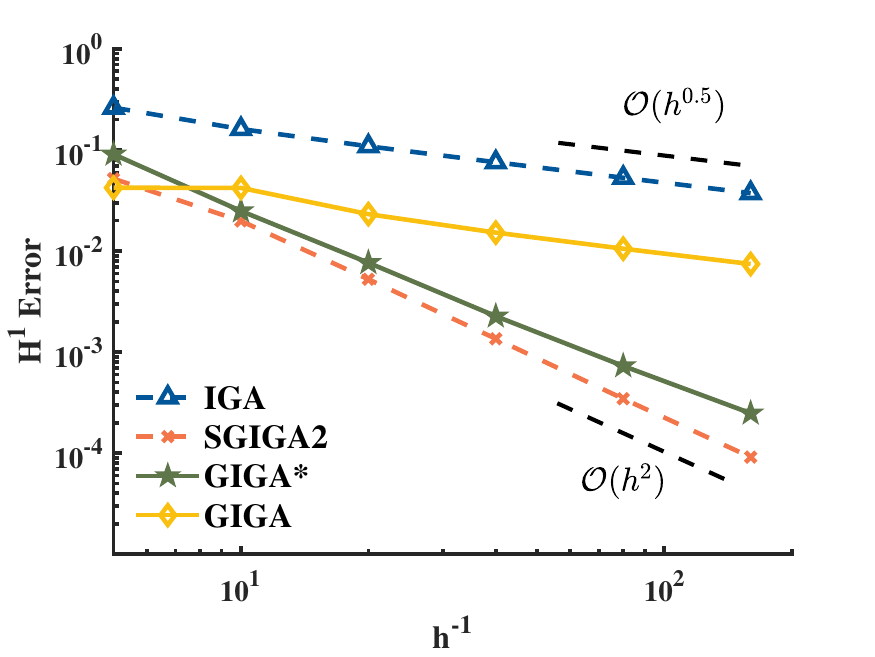}
        }\hfill
    \subfloat[$ a_0 = 100, a_1 = 1$]{
        \includegraphics[width=0.45\linewidth]{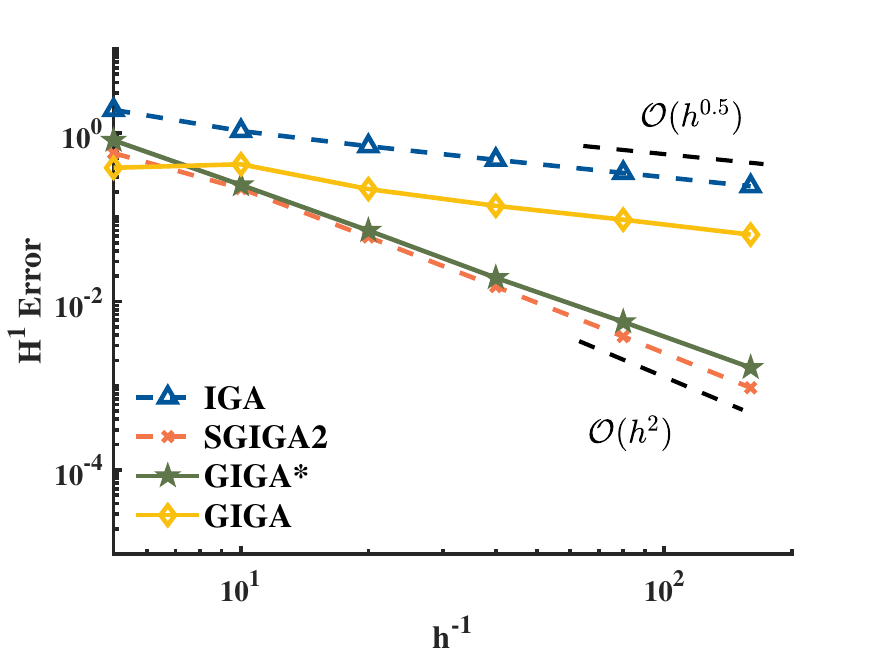}
        }
    \caption{$H^1$ errors with respect to $h$ (from $\frac{1}{5}$ to $\frac{1}{160}$) for an arc interface.}
    \label{the H1 error of yuanpan interface}
\end{figure}

\begin{figure}[htbp]
    \centering
    \subfloat[$a_0 = 20, a_1 = 1$]{     
        \includegraphics[width=0.45\linewidth]{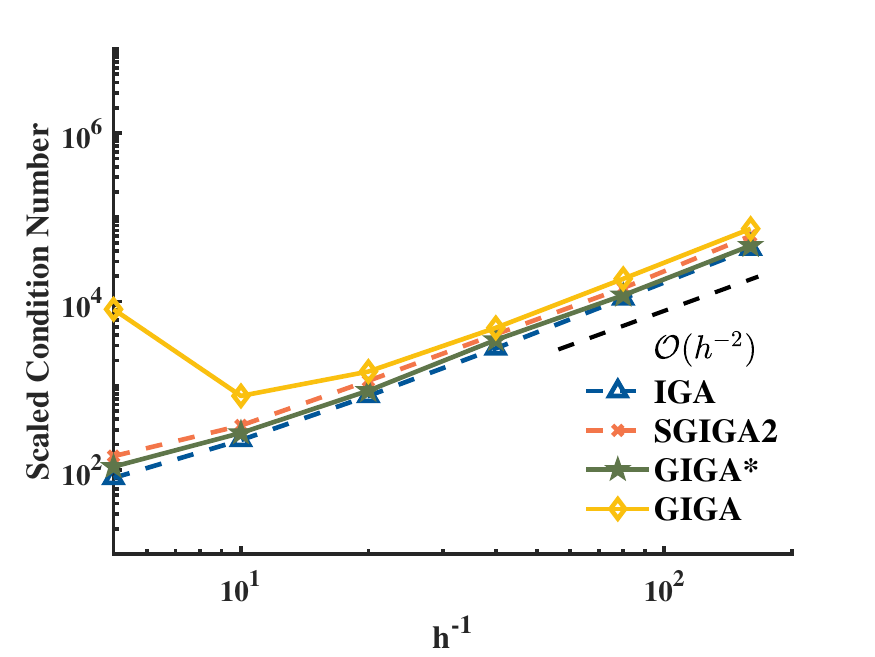}
        }\hfill
    \subfloat[$ a_0 = 100, a_1 = 1$]{
        \includegraphics[width=0.45\linewidth]{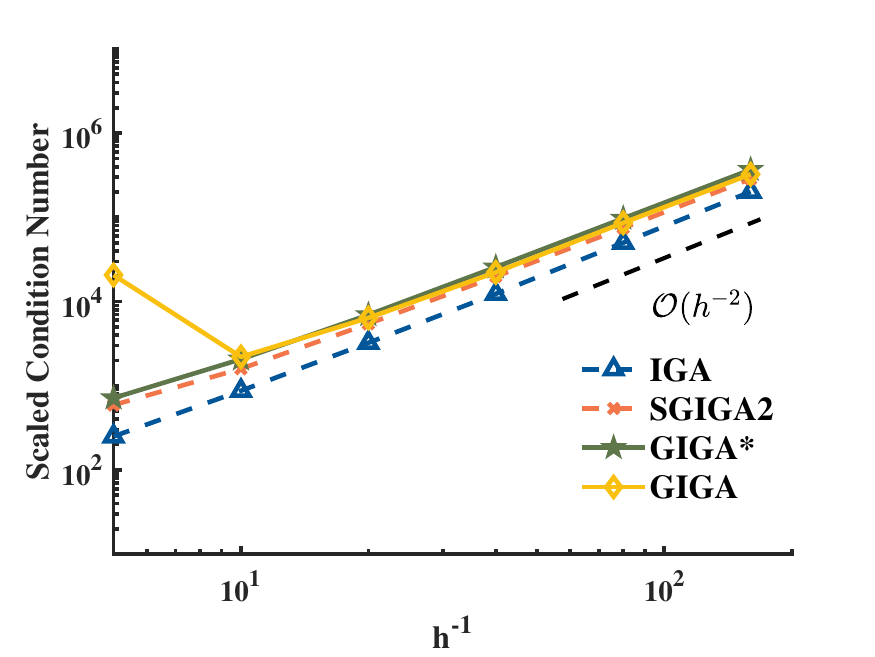}
        }
    \caption{The scaled condition numbers with respect to $h$ (from $ \frac{1}{5}$ to $\frac{1}{160}$) for an arc
interface.}
    \label{yuanpan_SCN}
\end{figure}

Fig~\ref{the L2 error of yuanpan interface}, Fig~\ref{the H1 error of yuanpan interface} and Fig~\ref{yuanpan_SCN} present the $ L^2 $ norm and $ H^1 $ seminorm errors and SCNs of several methods for an arc line interface. It is clearly observed that the SGIGA2 yields optimal convergence while the $ L^2 $ norm errors are $ O(h^3) $ and the $ H^1 $ seminorm errors are $ O(h^2) $. The SCNs of IGA, GIGA, GIGA* and SGIGA2 increase with the same order $ O(h^{-2}) $.

{The ultimate numerical investigation evaluates robustness of the proposed methodology through sensitivity analysis with respect to interface-boundary proximity.} We consider an example with a mobile interface $ y = \delta $, where $ \delta = 0.05\times 2^{-j}, j = 1,2,\dots,20$, and $ h = \frac{1}{20} $.
In this case, the interface $ \Gamma $ approaches the boundary of elements ($ \Omega = [0,1] \times [0,1] $) as $ \delta $ decreases. Letting $ c = \frac{a_0}{a_1} $ denote the contrast of the coefficients, we tested two cases, one with $ c = 10 $ and one with $ c = 100 $, setting $ a_1 = 1 $.

\begin{figure}[htbp]
    \centering
    \subfloat[$ a_0 = 10, a_1 = 1$]{     
        \includegraphics[width=0.45\linewidth]{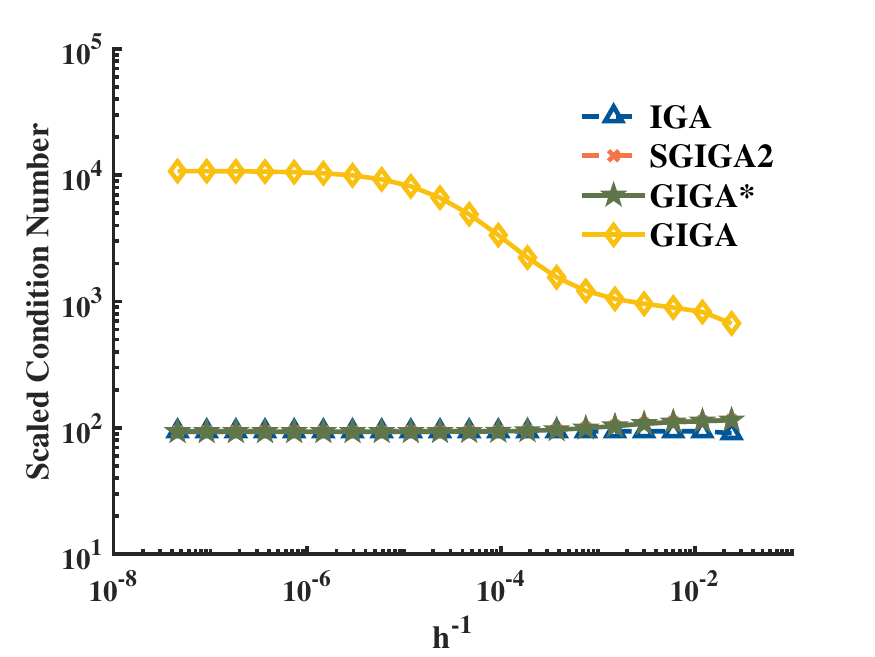}
        }\hfill
    \subfloat[$ a_0 = 100, a_1 = 1$]{
        \includegraphics[width=0.45\linewidth]{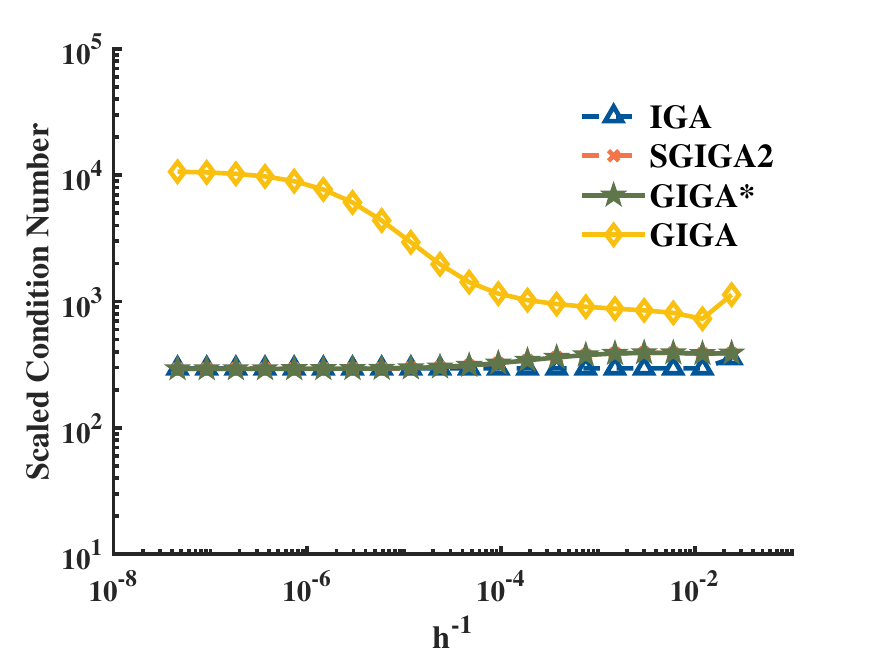}
        }
    \caption{Comparisons on the SCN with respect to $ \delta $, with fixed $ h = 1/20 $. The interface $ y = \delta $ approaches the boundary of element $ x = 0 $ as $ \delta $ decreases.}
    \label{robustness_test}
\end{figure}

Fig~\ref{robustness_test} shows that the SCNs of IGA, SGIGA2 and GIGA* do not grow as $ \delta $ decreases (the interface rapidly approaches the boundary of $ \Omega $ ) with fixed $ h = \frac{1}{20} $, which illustrates that they are robust. The SCN of GIGA increases
rapidly as $ \delta $ decreases, which means that GIGA is not robust.

\section{Conclusion and future work}\label{sec7}
{This paper presents the SGIGA2 method for solving smooth interface problems, achieving optimal convergence rates while eliminating the need for knot repetition. Through targeted enrichment space modifications, numerical experiments establish SGIGA2's stability and robustness.} The SCN of SGIGA2 does not depend on the location of the interface and increases with the order $O(h^{-2})$, which is the typical order of the SCN of IGA. Table~\ref{tab:Compare several method} {presents a comparative analysis of key methodologies.}
\begin{table*}[htbp]
    \centering
    \begin{tabular}{ccccccc}
   \toprule
   Prop. & FEM & IGA & \multicolumn{2}{c}{SGFEM2}  & GIGA* &SGIGA2   \\
   \cline{4-5}
   & & & \small{without LPCA} &\small{LPCA}& & \\
   \midrule

   $C^1$ continuity & $\times$ & $\checkmark$ & $\times$ & $\times$ & $\checkmark$ &$\checkmark$ \\
Optimal convergence & $\times$ & $\times$ &  $\checkmark$ & $\checkmark$& $\times$ &$\checkmark$\\
Stability & $\checkmark$ & $\checkmark$ &  $\times$ & $\checkmark$& $\checkmark$&$\checkmark$\\
Robustness& $\checkmark$ & $\checkmark$ & $\times$ & $\checkmark$& $\checkmark$&$\checkmark$\\
   \bottomrule
\end{tabular}
    \caption{A comparison of the properties of different methods.}
    \label{tab:Compare several method}
\end{table*}

{Although both \cite{ZHANG2020112889} and the present paper focus on biquadratic piecewise polynomials for interface problem approximation with all the desirable properties. The proposed methodology demonstrates several advantages over the SGFEM2 approach in~\cite{ZHANG2020112889}. Our formulation circumvents the need for local principal component analysis (LPCA) to mitigate stiffness matrix ill-conditioning, whereas SGFEM2 requires LPCA-based enrichment pruning through eigenvalue thresholding.} The LPCA method computes principal component coefficients
using the PCA method for each enriched node and then applies a percentage threshold to remove unnecessary principal component
enrichments. For example, in the circular interface example, the SGFEM2 removed more than 400 DOFs in the $160\times160$
mesh (Table~\ref{tab:Compare LPCA and our method}). Second, the present approach has less DOFs for the similar error. As shown
in Table~\ref{tab:Compare LPCA and our method}, we only need $30\%$ DOFs to have a similar error compared with SGFEM2.
\begin{table*}[htbp]
    \centering
    \begin{tabular}{ccccccc}
   \toprule
   N & FEM & IGA & \multicolumn{2}{c}{SGFEM2}  & GIGA* &SGIGA2   \\
   \cline{4-5}
   & & & \small{without LPCA} &\small{LPCA}& & \\
   \midrule
   5$\times$5 & 121 & 49 & 193 & 193 & 71 &115\\
10$\times$10 & 441 & 144 &  585 & 583& 203 &312\\
20$\times$20 & 1681 & 484 &  1981 & 1980& 592&808\\
40$\times$40& 6561 & 1764 & 7161 & 7150& 1972&2388\\
80$\times$80& 25921 & 6724 & 27145 & 27086&7140 & 7972\\
160$\times$160 & 103041 & 26244 & 105489 & 105036&27064 &28704\\
   \bottomrule
\end{tabular}
    \caption{The DOFs of the several space of the circle interface problem.}
    \label{tab:Compare LPCA and our method}
\end{table*}

{For future work, it is expected more improvements can be achieved.} One direction is the application to the higher order approximations and higher dimensional problems. It is quite challenging because even the well understood SGIGA can only handle specific interfaces, even if the interface can be accurately represented by NURBS. In addition, problems with two or more interfaces, possibly with highly heterogeneous coefficients in this region, still need further research. {Another challenging but interesting task is to explore} the numerical performance of SGIGA2 near the extraordinary points, while the current algorithm primarily targets the regular quadrilateral meshes.

\section*{Acknowledgements}
\vspace{-0.4cm} The authors were supported by the Strategic Priority Research Program of the Chinese Academy of
Sciences Grant No. XDB0640000, the Key Grant Project of the NSF of China (No.12494552), the NSF of China (No.12471360).


\bibliography{SGIGA2_xinli}


\begin{thebibliography}{55}
\ifx \bisbn   \undefined \def \bisbn  #1{ISBN #1}\fi
\ifx \binits  \undefined \def \binits#1{#1}\fi
\ifx \bauthor  \undefined \def \bauthor#1{#1}\fi
\ifx \batitle  \undefined \def \batitle#1{#1}\fi
\ifx \bjtitle  \undefined \def \bjtitle#1{#1}\fi
\ifx \bvolume  \undefined \def \bvolume#1{\textbf{#1}}\fi
\ifx \byear  \undefined \def \byear#1{#1}\fi
\ifx \bissue  \undefined \def \bissue#1{#1}\fi
\ifx \bfpage  \undefined \def \bfpage#1{#1}\fi
\ifx \blpage  \undefined \def \blpage #1{#1}\fi
\ifx \burl  \undefined \def \burl#1{\textsf{#1}}\fi
\ifx \doiurl  \undefined \def \doiurl#1{\url{https://doi.org/#1}}\fi
\ifx \betal  \undefined \def \betal{\textit{et al.}}\fi
\ifx \binstitute  \undefined \def \binstitute#1{#1}\fi
\ifx \binstitutionaled  \undefined \def \binstitutionaled#1{#1}\fi
\ifx \bctitle  \undefined \def \bctitle#1{#1}\fi
\ifx \beditor  \undefined \def \beditor#1{#1}\fi
\ifx \bpublisher  \undefined \def \bpublisher#1{#1}\fi
\ifx \bbtitle  \undefined \def \bbtitle#1{#1}\fi
\ifx \bedition  \undefined \def \bedition#1{#1}\fi
\ifx \bseriesno  \undefined \def \bseriesno#1{#1}\fi
\ifx \blocation  \undefined \def \blocation#1{#1}\fi
\ifx \bsertitle  \undefined \def \bsertitle#1{#1}\fi
\ifx \bsnm \undefined \def \bsnm#1{#1}\fi
\ifx \bsuffix \undefined \def \bsuffix#1{#1}\fi
\ifx \bparticle \undefined \def \bparticle#1{#1}\fi
\ifx \barticle \undefined \def \barticle#1{#1}\fi
\bibcommenthead
\ifx \bconfdate \undefined \def \bconfdate #1{#1}\fi
\ifx \botherref \undefined \def \botherref #1{#1}\fi
\ifx \url \undefined \def \url#1{\textsf{#1}}\fi
\ifx \bchapter \undefined \def \bchapter#1{#1}\fi
\ifx \bbook \undefined \def \bbook#1{#1}\fi
\ifx \bcomment \undefined \def \bcomment#1{#1}\fi
\ifx \oauthor \undefined \def \oauthor#1{#1}\fi
\ifx \citeauthoryear \undefined \def \citeauthoryear#1{#1}\fi
\ifx \endbibitem  \undefined \def \endbibitem {}\fi
\ifx \bconflocation  \undefined \def \bconflocation#1{#1}\fi
\ifx \arxivurl  \undefined \def \arxivurl#1{\textsf{#1}}\fi
\csname PreBibitemsHook\endcsname

\bibitem[\protect\citeauthoryear{Hughes et~al.}{2005}]{HUGHES20054135}
\begin{barticle}
\bauthor{\bsnm{Hughes}, \binits{T.J.R.}},
\bauthor{\bsnm{Cottrell}, \binits{J.A.}},
\bauthor{\bsnm{Bazilevs}, \binits{Y.}}:
\batitle{Isogeometric analysis: {CAD, finite elements, NURBS, exact geometry and mesh refinement}}.
\bjtitle{Computer Methods in Applied Mechanics and Engineering}
\bvolume{194}(\bissue{39}),
\bfpage{4135}--\blpage{4195}
(\byear{2005})
\end{barticle}
\endbibitem

\bibitem[\protect\citeauthoryear{Giannelli et~al.}{2012}]{GIANNELLI2012485}
\begin{barticle}
\bauthor{\bsnm{Giannelli}, \binits{C.}},
\bauthor{\bsnm{Juttler}, \binits{B.}},
\bauthor{\bsnm{Speleers}, \binits{H.}}:
\batitle{Thb-splines: The truncated basis for hierarchical splines}.
\bjtitle{Computer Aided Geometric Design}
\bvolume{29}(\bissue{7}),
\bfpage{485}--\blpage{498}
(\byear{2012})
\end{barticle}
\endbibitem

\bibitem[\protect\citeauthoryear{Wei et~al.}{2022}]{WEI2022114494}
\begin{barticle}
\bauthor{\bsnm{Wei}, \binits{X.}},
\bauthor{\bsnm{Li}, \binits{X.}},
\bauthor{\bsnm{Qian}, \binits{K.}},
\bauthor{\bsnm{Hughes}, \binits{T.J.R.}},
\bauthor{\bsnm{Zhang}, \binits{Y.J.}},
\bauthor{\bsnm{Casquero}, \binits{H.}}:
\batitle{Analysis-suitable unstructured {T-splines}: Multiple extraordinary points per face}.
\bjtitle{Computer Methods in Applied Mechanics and Engineering}
\bvolume{391},
\bfpage{114494}
(\byear{2022})
\end{barticle}
\endbibitem

\bibitem[\protect\citeauthoryear{Bazilevs et~al.}{2010}]{BAZILEVS2010229}
\begin{barticle}
\bauthor{\bsnm{Bazilevs}, \binits{Y.}},
\bauthor{\bsnm{Calo}, \binits{V.M.}},
\bauthor{\bsnm{Cottrell}, \binits{J.A.}},
\bauthor{\bsnm{Evans}, \binits{J.A.}},
\bauthor{\bsnm{Hughes}, \binits{T.J.R.}},
\bauthor{\bsnm{Lipton}, \binits{S.}},
\bauthor{\bsnm{Scott}, \binits{M.A.}},
\bauthor{\bsnm{Sederberg}, \binits{T.W.}}:
\batitle{Isogeometric analysis using {T}-splines}.
\bjtitle{Computer Methods in Applied Mechanics and Engineering}
\bvolume{199}(\bissue{5}),
\bfpage{229}--\blpage{263}
(\byear{2010}).
\bcomment{Computational Geometry and Analysis}
\end{barticle}
\endbibitem

\bibitem[\protect\citeauthoryear{Patrizi et~al.}{2020}]{PATRIZI2020113230}
\begin{barticle}
\bauthor{\bsnm{Patrizi}, \binits{F.}},
\bauthor{\bsnm{Manni}, \binits{C.}},
\bauthor{\bsnm{Pelosi}, \binits{F.}},
\bauthor{\bsnm{Speleers}, \binits{H.}}:
\batitle{Adaptive refinement with locally linearly independent lr b-splines: Theory and applications}.
\bjtitle{Computer Methods in Applied Mechanics and Engineering}
\bvolume{369},
\bfpage{113230}
(\byear{2020})
\end{barticle}
\endbibitem

\bibitem[\protect\citeauthoryear{Moes et~al.}{2003}]{MOES20033163}
\begin{barticle}
\bauthor{\bsnm{Moes}, \binits{N.}},
\bauthor{\bsnm{Cloirec}, \binits{M.}},
\bauthor{\bsnm{Cartraud}, \binits{P.}},
\bauthor{\bsnm{Remacle}, \binits{J.F.}}:
\batitle{A computational approach to handle complex microstructure geometries}.
\bjtitle{Computer Methods in Applied Mechanics and Engineering}
\bvolume{192}(\bissue{28}),
\bfpage{3163}--\blpage{3177}
(\byear{2003})
\end{barticle}
\endbibitem

\bibitem[\protect\citeauthoryear{Sauerland and Fries}{2011}]{SAUERLAND20113369}
\begin{barticle}
\bauthor{\bsnm{Sauerland}, \binits{H.}},
\bauthor{\bsnm{Fries}, \binits{T.-P.}}:
\batitle{The extended finite element method for two-phase and free-surface flows: A systematic study}.
\bjtitle{Journal of Computational Physics}
\bvolume{230}(\bissue{9}),
\bfpage{3369}--\blpage{3390}
(\byear{2011})
\end{barticle}
\endbibitem

\bibitem[\protect\citeauthoryear{Otmani et~al.}{2024}]{otmani2024accelerating}
\begin{botherref}
\oauthor{\bsnm{Otmani}, \binits{K.-E.}},
\oauthor{\bsnm{Mateo-Gab{\'\i}n}, \binits{A.}},
\oauthor{\bsnm{Rubio}, \binits{G.}},
\oauthor{\bsnm{Ferrer}, \binits{E.}}:
Accelerating high order discontinuous galerkin solvers through a clustering-based viscous/turbulent-inviscid domain decomposition.
Engineering with Computers,
1--16
(2024)
\end{botherref}
\endbibitem

\bibitem[\protect\citeauthoryear{R\v{e}thor\v{e} et~al.}{2007}]{RETHORE20075016}
\begin{barticle}
\bauthor{\bsnm{R\v{e}thor\v{e}}, \binits{J.}},
\bauthor{\bsnm{Hild}, \binits{F.}},
\bauthor{\bsnm{Roux}, \binits{S.}}:
\batitle{Shear-band capturing using a multiscale extended digital image correlation technique}.
\bjtitle{Computer Methods in Applied Mechanics and Engineering}
\bvolume{196}(\bissue{49}),
\bfpage{5016}--\blpage{5030}
(\byear{2007})
\end{barticle}
\endbibitem

\bibitem[\protect\citeauthoryear{Shakur}{2024}]{shakur2024isogeometric}
\begin{botherref}
\oauthor{\bsnm{Shakur}, \binits{E.}}:
Isogeometric analysis for solving discontinuous two-phase engineering problems with precise and explicit interface representation.
Engineering with Computers,
1--34
(2024)
\end{botherref}
\endbibitem

\bibitem[\protect\citeauthoryear{Gholampour et~al.}{2021}]{gholampour2021global}
\begin{barticle}
\bauthor{\bsnm{Gholampour}, \binits{F.}},
\bauthor{\bsnm{Hesameddini}, \binits{E.}},
\bauthor{\bsnm{Taleei}, \binits{A.}}:
\batitle{A global rbf-qr collocation technique for solving two-dimensional elliptic problems involving arbitrary interface}.
\bjtitle{Engineering with Computers}
\bvolume{37}(\bissue{4}),
\bfpage{3793}--\blpage{3811}
(\byear{2021})
\end{barticle}
\endbibitem

\bibitem[\protect\citeauthoryear{Hansbo and Hansbo}{2002}]{nitsche}
\begin{barticle}
\bauthor{\bsnm{Hansbo}, \binits{A.}},
\bauthor{\bsnm{Hansbo}, \binits{P.}}:
\batitle{An unfitted finite element method based on nitsches method for elliptic interface problems}.
\bjtitle{Computer Methods in Applied Mechanics and Engineering}
\bvolume{191(47-48)},
\bfpage{5537}--\blpage{5552}
(\byear{2002})
\end{barticle}
\endbibitem

\bibitem[\protect\citeauthoryear{Harari and Dolbow}{2010}]{harari2010analysis}
\begin{barticle}
\bauthor{\bsnm{Harari}, \binits{I.}},
\bauthor{\bsnm{Dolbow}, \binits{J.}}:
\batitle{Analysis of an efficient finite element method for embedded interface problems}.
\bjtitle{Computational Mechanics}
\bvolume{46},
\bfpage{205}--\blpage{211}
(\byear{2010})
\end{barticle}
\endbibitem

\bibitem[\protect\citeauthoryear{Huang et~al.}{2017}]{huang2017unfitted}
\begin{barticle}
\bauthor{\bsnm{Huang}, \binits{P.}},
\bauthor{\bsnm{Wu}, \binits{H.}},
\bauthor{\bsnm{Xiao}, \binits{Y.}}:
\batitle{An unfitted interface penalty finite element method for elliptic interface problems}.
\bjtitle{Computer Methods in Applied Mechanics and Engineering}
\bvolume{323},
\bfpage{439}--\blpage{460}
(\byear{2017})
\end{barticle}
\endbibitem

\bibitem[\protect\citeauthoryear{Chen and Liu}{2023}]{chen2023arbitrarily}
\begin{barticle}
\bauthor{\bsnm{Chen}, \binits{Z.}},
\bauthor{\bsnm{Liu}, \binits{Y.}}:
\batitle{An arbitrarily high order unfitted finite element method for elliptic interface problems with automatic mesh generation}.
\bjtitle{Journal of Computational Physics}
\bvolume{491},
\bfpage{112384}
(\byear{2023})
\end{barticle}
\endbibitem

\bibitem[\protect\citeauthoryear{Lehrenfeld and Reusken}{2013}]{lehrenfeld2013analysis}
\begin{barticle}
\bauthor{\bsnm{Lehrenfeld}, \binits{C.}},
\bauthor{\bsnm{Reusken}, \binits{A.}}:
\batitle{Analysis of a nitsche xfem-dg discretization for a class of two-phase mass transport problems}.
\bjtitle{SIAM journal on numerical analysis}
\bvolume{51}(\bissue{2}),
\bfpage{958}--\blpage{983}
(\byear{2013})
\end{barticle}
\endbibitem

\bibitem[\protect\citeauthoryear{Burman and Hansbo}{2010}]{burman2010fictitious}
\begin{barticle}
\bauthor{\bsnm{Burman}, \binits{E.}},
\bauthor{\bsnm{Hansbo}, \binits{P.}}:
\batitle{Fictitious domain finite element methods using cut elements: I. a stabilized lagrange multiplier method}.
\bjtitle{Computer Methods in Applied Mechanics and Engineering}
\bvolume{199}(\bissue{41-44}),
\bfpage{2680}--\blpage{2686}
(\byear{2010})
\end{barticle}
\endbibitem

\bibitem[\protect\citeauthoryear{Burman and Hansbo}{2012}]{burman2012fictitious}
\begin{barticle}
\bauthor{\bsnm{Burman}, \binits{E.}},
\bauthor{\bsnm{Hansbo}, \binits{P.}}:
\batitle{Fictitious domain finite element methods using cut elements: Ii. a stabilized nitsche method}.
\bjtitle{Applied Numerical Mathematics}
\bvolume{62}(\bissue{4}),
\bfpage{328}--\blpage{341}
(\byear{2012})
\end{barticle}
\endbibitem

\bibitem[\protect\citeauthoryear{Johansson and Larson}{2013}]{johansson2013high}
\begin{barticle}
\bauthor{\bsnm{Johansson}, \binits{A.}},
\bauthor{\bsnm{Larson}, \binits{M.G.}}:
\batitle{A high order discontinuous galerkin nitsche method for elliptic problems with fictitious boundary}.
\bjtitle{Numerische Mathematik}
\bvolume{123},
\bfpage{607}--\blpage{628}
(\byear{2013})
\end{barticle}
\endbibitem

\bibitem[\protect\citeauthoryear{Deng and Calo}{2020}]{DENG2020112558}
\begin{barticle}
\bauthor{\bsnm{Deng}, \binits{Q.}},
\bauthor{\bsnm{Calo}, \binits{V.}}:
\batitle{Higher order stable generalized finite element method for the elliptic eigenvalue and source problems with an interface in 1d}.
\bjtitle{Journal of Computational and Applied Mathematics}
\bvolume{368},
\bfpage{112558}
(\byear{2020})
\end{barticle}
\endbibitem

\bibitem[\protect\citeauthoryear{Babu\v{s}ka et~al.}{2017}]{BABUSKA201758}
\begin{barticle}
\bauthor{\bsnm{Babu\v{s}ka}, \binits{I.}},
\bauthor{\bsnm{Banerjee}, \binits{U.}},
\bauthor{\bsnm{Kergrene}, \binits{K.}}:
\batitle{Strongly stable generalized finite element method: Application to interface problems}.
\bjtitle{Computer Methods in Applied Mechanics and Engineering}
\bvolume{327},
\bfpage{58}--\blpage{92}
(\byear{2017})
\end{barticle}
\endbibitem

\bibitem[\protect\citeauthoryear{Zhang et~al.}{2019}]{ZHANG2019538}
\begin{barticle}
\bauthor{\bsnm{Zhang}, \binits{Q.}},
\bauthor{\bsnm{Banerjee}, \binits{U.}},
\bauthor{\bsnm{Babu\v{s}ka}, \binits{I.}}:
\batitle{Strongly stable generalized finite element method {(SSGFEM)} for a non-smooth interface problem}.
\bjtitle{Computer Methods in Applied Mechanics and Engineering}
\bvolume{344},
\bfpage{538}--\blpage{568}
(\byear{2019})
\end{barticle}
\endbibitem

\bibitem[\protect\citeauthoryear{Zhang et~al.}{2020}]{ZHANG2020112926}
\begin{barticle}
\bauthor{\bsnm{Zhang}, \binits{Q.}},
\bauthor{\bsnm{Banerjee}, \binits{U.}},
\bauthor{\bsnm{Babu\v{s}ka}, \binits{I.}}:
\batitle{Strongly stable generalized finite element method {(SSGFEM)} for a non-smooth interface problem {II}: A simplified algorithm}.
\bjtitle{Computer Methods in Applied Mechanics and Engineering}
\bvolume{363},
\bfpage{112926}
(\byear{2020})
\end{barticle}
\endbibitem

\bibitem[\protect\citeauthoryear{Zhu et~al.}{2020}]{ZHU2020112475}
\begin{barticle}
\bauthor{\bsnm{Zhu}, \binits{P.}},
\bauthor{\bsnm{Zhang}, \binits{Q.}},
\bauthor{\bsnm{Liu}, \binits{T.}}:
\batitle{Stable generalized finite element method {(SGFEM)} for parabolic interface problems}.
\bjtitle{Journal of Computational and Applied Mathematics}
\bvolume{367},
\bfpage{112475}
(\byear{2020})
\end{barticle}
\endbibitem

\bibitem[\protect\citeauthoryear{Gong et~al.}{2024}]{GONG2024115540}
\begin{barticle}
\bauthor{\bsnm{Gong}, \binits{W.}},
\bauthor{\bsnm{Li}, \binits{H.}},
\bauthor{\bsnm{Zhang}, \binits{Q.}}:
\batitle{Improved enrichments and numerical integrations in {SGFEM} for interface problems}.
\bjtitle{Journal of Computational and Applied Mathematics}
\bvolume{438},
\bfpage{115540}
(\byear{2024})
\end{barticle}
\endbibitem

\bibitem[\protect\citeauthoryear{Jia et~al.}{2015}]{8172107}
\begin{barticle}
\bauthor{\bsnm{Jia}, \binits{Y.}},
\bauthor{\bsnm{Anitescu}, \binits{C.}},
\bauthor{\bsnm{Ghorashi}, \binits{S.S.}},
\bauthor{\bsnm{Rabczuk}, \binits{T.}}:
\batitle{Extended isogeometric analysis for material interface problems}.
\bjtitle{IMA Journal of Applied Mathematics}
\bvolume{80}(\bissue{3}),
\bfpage{608}--\blpage{633}
(\byear{2015})
\end{barticle}
\endbibitem

\bibitem[\protect\citeauthoryear{Zhang et~al.}{2022}]{ZHANG2022114053}
\begin{barticle}
\bauthor{\bsnm{Zhang}, \binits{J.}},
\bauthor{\bsnm{Deng}, \binits{Q.}},
\bauthor{\bsnm{Li}, \binits{X.}}:
\batitle{A generalized isogeometric analysis of elliptic eigenvalue and source problems with an interface}.
\bjtitle{Journal of Computational and Applied Mathematics}
\bvolume{407},
\bfpage{114053}
(\byear{2022})
\end{barticle}
\endbibitem

\bibitem[\protect\citeauthoryear{Tan et~al.}{2015}]{TAN20151382}
\begin{barticle}
\bauthor{\bsnm{Tan}, \binits{M.H.Y.}},
\bauthor{\bsnm{Safdari}, \binits{M.}},
\bauthor{\bsnm{Najafi}, \binits{A.R.}},
\bauthor{\bsnm{Geubelle}, \binits{P.H.}}:
\batitle{A {NURBS-based} interface-enriched generalized finite element scheme for the thermal analysis and design of microvascular composites}.
\bjtitle{Computer Methods in Applied Mechanics and Engineering}
\bvolume{283},
\bfpage{1382}--\blpage{1400}
(\byear{2015})
\end{barticle}
\endbibitem

\bibitem[\protect\citeauthoryear{Kiendl et~al.}{2009}]{KIENDL20093902}
\begin{barticle}
\bauthor{\bsnm{Kiendl}, \binits{J.}},
\bauthor{\bsnm{Bletzinger}, \binits{K.-U.}},
\bauthor{\bsnm{Linhard}, \binits{J.}},
\bauthor{\bsnm{Wuchner}, \binits{R.}}:
\batitle{Isogeometric shell analysis with {Kirchhoff-Love} elements}.
\bjtitle{Computer Methods in Applied Mechanics and Engineering}
\bvolume{198}(\bissue{49}),
\bfpage{3902}--\blpage{3914}
(\byear{2009})
\end{barticle}
\endbibitem

\bibitem[\protect\citeauthoryear{G\v{o}mez et~al.}{2008}]{GOMEZ20084333}
\begin{barticle}
\bauthor{\bsnm{G\v{o}mez}, \binits{H.}},
\bauthor{\bsnm{Calo}, \binits{V.M.}},
\bauthor{\bsnm{Bazilevs}, \binits{Y.}},
\bauthor{\bsnm{Hughes}, \binits{T.J.R.}}:
\batitle{Isogeometric analysis of the {Cahn-Hilliard} phase-field model}.
\bjtitle{Computer Methods in Applied Mechanics and Engineering}
\bvolume{197}(\bissue{49}),
\bfpage{4333}--\blpage{4352}
(\byear{2008})
\end{barticle}
\endbibitem

\bibitem[\protect\citeauthoryear{Cazzani et~al.}{2016}]{Cazzani2016IsogeometricAO}
\begin{barticle}
\bauthor{\bsnm{Cazzani}, \binits{A.}},
\bauthor{\bsnm{Malag{\`u}}, \binits{M.}},
\bauthor{\bsnm{Turco}, \binits{E.}}:
\batitle{Isogeometric analysis of plane-curved beams}.
\bjtitle{Mathematics and Mechanics of Solids}
\bvolume{21},
\bfpage{562}--\blpage{577}
(\byear{2016})
\end{barticle}
\endbibitem

\bibitem[\protect\citeauthoryear{Thai et~al.}{2013}]{THAI2013196}
\begin{barticle}
\bauthor{\bsnm{Thai}, \binits{C.H.}},
\bauthor{\bsnm{Ferreira}, \binits{A.J.M.}},
\bauthor{\bsnm{Carrera}, \binits{E.}},
\bauthor{\bsnm{Nguyen-Xuan}, \binits{H.}}:
\batitle{Isogeometric analysis of laminated composite and sandwich plates using a layerwise deformation theory}.
\bjtitle{Composite Structures}
\bvolume{104},
\bfpage{196}--\blpage{214}
(\byear{2013})
\end{barticle}
\endbibitem

\bibitem[\protect\citeauthoryear{Babu\v{s}ka et~al.}{2003}]{BabuOsborn_2003}
\begin{barticle}
\bauthor{\bsnm{Babu\v{s}ka}, \binits{I.}},
\bauthor{\bsnm{Banerjee}, \binits{U.}},
\bauthor{\bsnm{Osborn}, \binits{J.E.}}:
\batitle{Survey of meshless and generalized finite element methods: A unified approach}.
\bjtitle{Acta Numerica}
\bvolume{12},
\bfpage{1}--\blpage{125}
(\byear{2003})
\end{barticle}
\endbibitem

\bibitem[\protect\citeauthoryear{Fries and Belytschko}{2010}]{fries2010extended}
\begin{barticle}
\bauthor{\bsnm{Fries}, \binits{T.-P.}},
\bauthor{\bsnm{Belytschko}, \binits{T.}}:
\batitle{The extended/generalized finite element method: an overview of the method and its applications}.
\bjtitle{International journal for numerical methods in engineering}
\bvolume{84}(\bissue{3}),
\bfpage{253}--\blpage{304}
(\byear{2010})
\end{barticle}
\endbibitem

\bibitem[\protect\citeauthoryear{Abdelaziz and Hamouine}{2008}]{ABDELAZIZ20081141}
\begin{barticle}
\bauthor{\bsnm{Abdelaziz}, \binits{Y.}},
\bauthor{\bsnm{Hamouine}, \binits{A.}}:
\batitle{A survey of the extended finite element}.
\bjtitle{Computers Structures}
\bvolume{86}(\bissue{11}),
\bfpage{1141}--\blpage{1151}
(\byear{2008})
\end{barticle}
\endbibitem

\bibitem[\protect\citeauthoryear{Iqbal et~al.}{2022}]{iqbal2022enriched}
\begin{botherref}
\oauthor{\bsnm{Iqbal}, \binits{M.}},
\oauthor{\bsnm{Alam}, \binits{K.}},
\oauthor{\bsnm{Ahmad}, \binits{A.}},
\oauthor{\bsnm{Maqsood}, \binits{S.}},
\oauthor{\bsnm{Ullah}, \binits{H.}},
\oauthor{\bsnm{Ullah}, \binits{B.}}:
An enriched finite element method for efficient solutions of transient heat diffusion problems with multiple heat sources.
Engineering with Computers,
1--17
(2022)
\end{botherref}
\endbibitem

\bibitem[\protect\citeauthoryear{Babu\v{s}ka and Banerjee}{2012}]{BABUSKA201291}
\begin{barticle}
\bauthor{\bsnm{Babu\v{s}ka}, \binits{I.}},
\bauthor{\bsnm{Banerjee}, \binits{U.}}:
\batitle{Stable generalized finite element method {(SGFEM)}}.
\bjtitle{Computer Methods in Applied Mechanics and Engineering}
\bvolume{201-204},
\bfpage{91}--\blpage{111}
(\byear{2012})
\end{barticle}
\endbibitem

\bibitem[\protect\citeauthoryear{Agathos et~al.}{2019}]{AGATHOS20191051}
\begin{barticle}
\bauthor{\bsnm{Agathos}, \binits{K.}},
\bauthor{\bsnm{Bordas}, \binits{S.P.A.}},
\bauthor{\bsnm{Chatzi}, \binits{E.}}:
\batitle{Improving the conditioning of {XFEM/GFEM} for fracture mechanics problems through enrichment quasi-orthogonalization}.
\bjtitle{Computer Methods in Applied Mechanics and Engineering}
\bvolume{346},
\bfpage{1051}--\blpage{1073}
(\byear{2019})
\end{barticle}
\endbibitem

\bibitem[\protect\citeauthoryear{Schweitzer}{2011}]{schweitzer2011stable}
\begin{barticle}
\bauthor{\bsnm{Schweitzer}, \binits{M.A.}}:
\batitle{Stable enrichment and local preconditioning in the particle-partition of unity method}.
\bjtitle{Numerische Mathematik}
\bvolume{118}(\bissue{1}),
\bfpage{137}--\blpage{170}
(\byear{2011})
\end{barticle}
\endbibitem

\bibitem[\protect\citeauthoryear{Babu\v{s}ka et~al.}{1994}]{babu94}
\begin{barticle}
\bauthor{\bsnm{Babu\v{s}ka}, \binits{I.}},
\bauthor{\bsnm{Caloz}, \binits{G.}},
\bauthor{\bsnm{Osborn}, \binits{J.E.}}:
\batitle{Special finite element methods for a class of second order elliptic problems with rough coefficients}.
\bjtitle{SIAM Journal on Numerical Analysis}
\bvolume{31}(\bissue{4}),
\bfpage{945}--\blpage{981}
(\byear{1994})
\end{barticle}
\endbibitem

\bibitem[\protect\citeauthoryear{Sauerland and Fries}{2013}]{SAUERLAND201341}
\begin{barticle}
\bauthor{\bsnm{Sauerland}, \binits{H.}},
\bauthor{\bsnm{Fries}, \binits{T.-P.}}:
\batitle{The stable {XFEM} for two-phase flows}.
\bjtitle{Computers Fluids}
\bvolume{87},
\bfpage{41}--\blpage{49}
(\byear{2013})
\end{barticle}
\endbibitem

\bibitem[\protect\citeauthoryear{Gupta et~al.}{2013}]{GUPTA201323}
\begin{barticle}
\bauthor{\bsnm{Gupta}, \binits{V.}},
\bauthor{\bsnm{Duarte}, \binits{C.A.}},
\bauthor{\bsnm{Babu\v{s}ka}, \binits{I.}},
\bauthor{\bsnm{Banerjee}, \binits{U.}}:
\batitle{A stable and optimally convergent generalized {FEM (SGFEM)} for linear elastic fracture mechanics}.
\bjtitle{Computer Methods in Applied Mechanics and Engineering}
\bvolume{266},
\bfpage{23}--\blpage{39}
(\byear{2013})
\end{barticle}
\endbibitem

\bibitem[\protect\citeauthoryear{Zhang et~al.}{2016}]{ZHANG2016476}
\begin{barticle}
\bauthor{\bsnm{Zhang}, \binits{Q.}},
\bauthor{\bsnm{Babu\v{s}ka}, \binits{I.}},
\bauthor{\bsnm{Banerjee}, \binits{U.}}:
\batitle{Robustness in stable generalized finite element methods {(SGFEM)} applied to poisson problems with crack singularities}.
\bjtitle{Computer Methods in Applied Mechanics and Engineering}
\bvolume{311},
\bfpage{476}--\blpage{502}
(\byear{2016})
\end{barticle}
\endbibitem

\bibitem[\protect\citeauthoryear{Fries}{2008}]{article}
\begin{barticle}
\bauthor{\bsnm{Fries}, \binits{T.-P.}}:
\batitle{A corrected {XFEM} approximation without problems in blending elements}.
\bjtitle{International Journal for Numerical Methods in Engineering}
\bvolume{75},
\bfpage{503}--\blpage{532}
(\byear{2008})
\end{barticle}
\endbibitem

\bibitem[\protect\citeauthoryear{Cheng and Fries}{2009}]{Cheng2009HigherorderXF}
\begin{botherref}
\oauthor{\bsnm{Cheng}, \binits{K.W.}},
\oauthor{\bsnm{Fries}, \binits{T.P.}}:
Higher-order {XFEM} for curved strong and weak discontinuities.
International Journal for Numerical Methods in Engineering
\textbf{82}
(2009)
\end{botherref}
\endbibitem

\bibitem[\protect\citeauthoryear{Hou et~al.}{2020}]{HOU2020113135}
\begin{barticle}
\bauthor{\bsnm{Hou}, \binits{W.}},
\bauthor{\bsnm{Jiang}, \binits{K.}},
\bauthor{\bsnm{Zhu}, \binits{X.}},
\bauthor{\bsnm{Shen}, \binits{Y.}},
\bauthor{\bsnm{Li}, \binits{Y.}},
\bauthor{\bsnm{Zhang}, \binits{X.}},
\bauthor{\bsnm{Hu}, \binits{P.}}:
\batitle{Extended isogeometric analysis with strong imposing essential boundary conditions for weak discontinuous problems using {B++ splines}}.
\bjtitle{Computer Methods in Applied Mechanics and Engineering}
\bvolume{370},
\bfpage{113135}
(\byear{2020})
\end{barticle}
\endbibitem

\bibitem[\protect\citeauthoryear{{Durga Rao} and Raju}{2019}]{DURGARAO2019535}
\begin{barticle}
\bauthor{\bsnm{{Durga Rao}}, \binits{S.S.}},
\bauthor{\bsnm{Raju}, \binits{S.}}:
\batitle{{Stable Generalized Iso-Geometric Analysis (SGIGA)} for problems with discontinuities and singularities}.
\bjtitle{Computer Methods in Applied Mechanics and Engineering}
\bvolume{348},
\bfpage{535}--\blpage{574}
(\byear{2019})
\end{barticle}
\endbibitem

\bibitem[\protect\citeauthoryear{Hu et~al.}{2024}]{HU2024115792}
\begin{barticle}
\bauthor{\bsnm{Hu}, \binits{W.}},
\bauthor{\bsnm{Zhang}, \binits{J.}},
\bauthor{\bsnm{Li}, \binits{X.}}:
\batitle{Higher order stable generalized isogeometric analysis for interface problems}.
\bjtitle{Journal of Computational and Applied Mathematics}
\bvolume{444},
\bfpage{115792}
(\byear{2024})
\end{barticle}
\endbibitem

\bibitem[\protect\citeauthoryear{Zhang and Babu\v{s}ka}{2020}]{ZHANG2020112889}
\begin{barticle}
\bauthor{\bsnm{Zhang}, \binits{Q.}},
\bauthor{\bsnm{Babu\v{s}ka}, \binits{I.}}:
\batitle{A stable generalized finite element method {(SGFEM)} of degree two for interface problems}.
\bjtitle{Computer Methods in Applied Mechanics and Engineering}
\bvolume{363},
\bfpage{112889}
(\byear{2020})
\end{barticle}
\endbibitem

\bibitem[\protect\citeauthoryear{Kang et~al.}{2022}]{xinli_quasi}
\begin{botherref}
\oauthor{\bsnm{Kang}, \binits{H.}},
\oauthor{\bsnm{Yong}, \binits{Z.}},
\oauthor{\bsnm{Li}, \binits{X.}}:
Quasi-interpolation for analysis-suitable {T}-splines.
Computer aided geometric design
\textbf{98}
(2022)
\end{botherref}
\endbibitem

\bibitem[\protect\citeauthoryear{Barrera et~al.}{2020}]{barrera2020non}
\begin{barticle}
\bauthor{\bsnm{Barrera}, \binits{D.}},
\bauthor{\bsnm{El~Mokhtari}, \binits{F.}},
\bauthor{\bsnm{Ib{\'a}{\~n}ez}, \binits{M.J.}},
\bauthor{\bsnm{Sbibih}, \binits{D.}}:
\batitle{Non-uniform quasi-interpolation for solving {Hammerstein} integral equations}.
\bjtitle{International Journal of Computer Mathematics}
\bvolume{97}(\bissue{1-2}),
\bfpage{72}--\blpage{84}
(\byear{2020})
\end{barticle}
\endbibitem

\bibitem[\protect\citeauthoryear{Haasemann et~al.}{2011}]{haasemann2011development}
\begin{barticle}
\bauthor{\bsnm{Haasemann}, \binits{G.}},
\bauthor{\bsnm{K{\"a}stner}, \binits{M.}},
\bauthor{\bsnm{Pr{\"u}ger}, \binits{S.}},
\bauthor{\bsnm{Ulbricht}, \binits{V.}}:
\batitle{Development of a quadratic finite element formulation based on the {XFEM and NURBS}}.
\bjtitle{International Journal for Numerical Methods in Engineering}
\bvolume{86}(\bissue{4-5}),
\bfpage{598}--\blpage{617}
(\byear{2011})
\end{barticle}
\endbibitem

\bibitem[\protect\citeauthoryear{K{\"a}stner et~al.}{2013}]{kastner2013higher}
\begin{barticle}
\bauthor{\bsnm{K{\"a}stner}, \binits{M.}},
\bauthor{\bsnm{M{\"u}ller}, \binits{S.}},
\bauthor{\bsnm{Goldmann}, \binits{J.}},
\bauthor{\bsnm{Spieler}, \binits{C.}},
\bauthor{\bsnm{Brummund}, \binits{J.}},
\bauthor{\bsnm{Ulbricht}, \binits{V.}}:
\batitle{Higher-order extended {FEM} for weak discontinuities level set representation, quadrature and application to magneto-mechanical problems}.
\bjtitle{International Journal for Numerical Methods in Engineering}
\bvolume{93}(\bissue{13}),
\bfpage{1403}--\blpage{1424}
(\byear{2013})
\end{barticle}
\endbibitem

\bibitem[\protect\citeauthoryear{Fries and Belytschko}{2010}]{nme2914}
\begin{barticle}
\bauthor{\bsnm{Fries}, \binits{T.-P.}},
\bauthor{\bsnm{Belytschko}, \binits{T.}}:
\batitle{The extended/generalized finite element method: An overview of the method and its applications}.
\bjtitle{International Journal for Numerical Methods in Engineering}
\bvolume{84}(\bissue{3}),
\bfpage{253}--\blpage{304}
(\byear{2010})
\end{barticle}
\endbibitem

\bibitem[\protect\citeauthoryear{Fries}{2008}]{fries2008corrected}
\begin{barticle}
\bauthor{\bsnm{Fries}, \binits{T.-P.}}:
\batitle{A corrected {XFEM} approximation without problems in blending elements}.
\bjtitle{International Journal for Numerical Methods in Engineering}
\bvolume{75}(\bissue{5}),
\bfpage{503}--\blpage{532}
(\byear{2008})
\end{barticle}
\endbibitem

\end{thebibliography}

\end{document}